\documentclass[12pt]{article}
\usepackage{amscd,amssymb,epic,amsmath}
\usepackage[mathscr]{eucal}
\usepackage[all]{xy}
\usepackage{theorem,marginnote}

\newtheorem{proposition}{Proposition}
\newtheorem{lemma}{Lemma}
\newtheorem{corollary}{Corollary}
\newtheorem{theorem}{Theorem}

\newtheorem{remark}{Remark}

\begin{document}

\title {Noncommutative Harmonic Analysis on Quantum Hyperbolic Spaces. The Laplace-Beltrami Operator} 
\author{Olga Bershteyn}

\date{Institute of Mathematical Studies, University of Copenhagen \\
e-mail: olga@math.ku.dk, olya.bersht@gmail.com}



\maketitle

\begin{abstract}
In this paper we study the Laplace-Beltrami operator on quantum complex hyperbolic spaces. We describe its action in terms of certain $q$-difference operators of second order and prove spectral theorems for these operators. The corresponding eigenfunctions are related to Al-Salam-Chihara polynomials. The obtained Plancherel measure is related to a quantum analog for the Harish-Chandra c-function.
\end{abstract}

\section{Introduction}

Harmonic analysis on real reductive groups and their homogeneous spaces attracted a lot of prominent mathematicians during the last 50 years. One of the reasons is its deep connections with representation theory, theory of special functions and differential equations. This knot of problems and links is so atractive that during the last twenty years several independent groups of mathematicians successfully elaborated similar theories in other branches of mathematics. That happened within the framework of quantum group theory as well. One could mention the papers of Ueno \cite{U}, Vaksman and Korogodskii \cite{KV}, Koelink and Stokman \cite{KS}, and Ip \cite{I}, where some types of quantum groups and quantum homogeneous spaces were considered. In this paper we contribute to the list of elaborated examples by considering so called quantum complex hyperbolic spaces. 

Let us recall the classical results on complex hyperbolic spaces. Let $m,n\in\mathbb{N}$, $m\ge 2$, and $N\overset{\mathrm{def}}{=}n+m$. Let us consider the hyperboloid  
$$
\widehat{\mathscr{H}}_{n,m}=\left\{(t_1,\ldots,t_N)\in\mathbb{C}^N\left|\:
-\sum^n_{j=1}|t_j|^2+\sum^N_{j=n+1}|t_{j}|^2=1\right.\right\}.
$$
One can identify points of $\widehat{\mathscr{H}}_{n,m}$ by the equivalence relation $a \sim b$ if $a=e^{i\theta} b$ for some $\theta \in [0,2\pi]$. It is well-known that $\mathscr{H}_{n,m}=\widehat{\mathscr{H}}_{n,m}/\sim$ is a homogeneous space of the real reductive group $SU_{n,m}$, namely
$\mathscr{H}_{n,m}=SU_{n,m}/S(U_{n,m-1}\times U_1)$. It is called a complex hyperbolic space and is one of well studied and important examples of pseudo-Hermitian symmetric spaces of rank 1. The Faraut paper \cite{Faraut} on hyperbolic spaces made a great impact into the theory of
semisimple symmetric spaces of rank 1. In particular, Faraut obtained the Plancherel type decomposition for the quasi-regular representation of $SU_{n,m}$ in $L^2(\mathscr{H}_{n,m}, d\nu)$ in unitary irreducibles ($d\nu$ is an invariant measure on $\mathscr{H}_{n,m}$). The main his tools were so called spherical distributions and decomposition of the Laplace-Beltrami operator. Also there are numerous papers of Molchanov, van Dijk, Schlichtkrull and others (see \cite{Molch1,Molch2,Dijk-Sh,Schl} and references therein) on representation theory related to complex hyperbolic spaces and harmonic analysis on them.

In the last decade of 20th century Leonid Vaksman with colleagues started to study harmonic analysis on quantum symmetric spaces. It appeared back then that several problems of classical complex and harmonic analysis on symmetric spaces have quantum analogs. They focused mostly on Hermitian symmetric spaces of non compact type. The idea was to obtain a Plancherel type decomposition for quasi-regular representation of a quantum universal enveloping algebra $U_q \mathfrak{g}$, and the main tool was to study decompositions of certain invariant operators. Unfortunately, the quantum case was developed in full only for the simplest case of the group $SU(1,1)$ and its quantum counterpart $U_q \mathfrak{su}_{1,1}$. The corresponding quantum symmetric space is called the quantum unit disc and it is the simplest Hermitian symmetric space. For this special case, the Laplacian was constructed, its eigenfunctions were obtained, and this gives then the spectral decomposition for the Laplacian and for the whole quasi-regular representaion in the quantum unit disc (See \cite{V} for the summary and references therein). Later some intermediate results were obtained for other Hermitian symmetric spaces. In their paper \cite{SZ}, Shklyarov and Zhang obtained spectral decomposition for the Laplace operator on quantum unit balls in $\mathbb C^n$. Certain similar problems were solved for the case of quantum matrix balls which are quantum analogs of the symmetric spaces $SU(n,n)/S(U_n \times U_n)$ (see \cite{BK}). In this case it was possible to construct a family of commuting invariant operators and obtain their spectral decompositions. The present paper and the previous one \cite{BK1} containes the same kind of results for another class of symmetric spaces.

The paper is organised as follows. Section \ref{sec_1} contains background from the general theory of quantum groups. The key objects there are the Hopf dual algebras $U_q \mathfrak{sl}_N$ and $\mathbb C[SL_N]_q$. Important additional structures are the noncompact involutions on both these algebras. The algebra $U_q \mathfrak{sl}_N$ together with the corresponding noncompact involution is treated as a quantum analogs for the real group $SU(n,m)$. In Section \ref{sec_2} we recall the basic notions and notation regarding quantum complex hyperbolic spaces, following \cite{BS} and \cite{BK1}. Namely, we introduce a quantum analog $\mathcal{D}(\mathscr{H}_{n,m})_{q,\mathfrak{k}}$ for the algebra of ($K$-finite, where $K=S(U_n \times U_m)$) smooth compactly supported functions on $\mathscr{H}_{n,m}$ and provide an invariant integral (i.e., a covariant linear functional) on it. We give a description of $U_q \mathfrak{k}$-isotypic components of  $\mathcal{D}(\mathscr{H}_{n,m})_{q,\mathfrak{k}}$ in Proposition \ref{Kdecomp_H}.  The Laplace-Beltrami operator $\Box_q$ is constructed in Section \ref{sec_6}. Furthermore we reduce the action of $\Box_q$ to the $U_q \mathfrak{k}$- components of $\mathcal{D}(\mathscr{H}_{n,m})_{q,\mathfrak{k}}$. This action can be essentially described by means of certain $q$-difference operators $A^{\Lambda,\Lambda'}$ of second order. Explicit formula for these operators is given in Proposition \ref{Laplace_on_highestvector}. In the next section we study these bounded self-adjoint operators $A^{\Lambda,\Lambda'}$ and find their eigenfunctions. A remarkable property of these eigenfunctions is that they can be written as Al-Salam-Chihara polynomials (in a dual variable), so the spectral analysis of $A^{\Lambda,\Lambda'}$ is reduced to some well-known facts about the $q$-special functions. Section \ref{sec_9} contains spectral theorems for the $q$-difference operators $A^{\Lambda,\Lambda'}$. We also obtain the Plancherel measure, and describe its connection with the Harish-Chandra $c$-function. A side remark is that the Plancherel measure gives us a clue about which unitary representaions should appear in the decomposition of the left quasiregular one. Appendix contains proofs of the most technical statements.

The author would like to express her gratitude to L.~Vaksman and D.~Shklyarov for sharing their ideas and first results. The author is also grateful to H.P.~Jakobsen for many valuable discussions.

\section{Preliminaries on quantum group theory}\label{sec_1}

Everywhere in the sequel we suppose $q\in(0,1)$. All algebras are associative and unital. Let us recall some standard notions and notations from the theory of quantum groups. All the facts below can be found in standard books on quantum group theory \cite{Jant,Kl-Sch}.

The Hopf algebra $U_q\mathfrak{sl}_{N}$ is given as an associative algebra with generators $K_i$, $K_i^{-1}$, $E_i$, $F_i$, $i=1,2,\ldots,N-1$, and the well-known Drinfeld-Jimbo relations (see \cite{Jant})
$$ K_iK_j=K_jK_i,\qquad K_iK_i^{-1}=K_i^{-1}K_i=1,$$
$$K_iE_j=q^{a_{ij}}E_jK_i, \quad K_iF_j=q^{-a_{ij}}F_jK_i, \qquad a_{ij}=\begin{cases} 2, & i=j \\ -1, & |i-j|=1, \\ 0, & \text{otherwise}. \end{cases}$$
$$E_iF_j-F_jE_i=\delta_{ij}(K_i-K_i^{-1})/(q-q^{-1}), $$
\begin{equation*}
E_i^2E_j-(q+q^{-1})E_iE_jE_i+E_jE_i^2=0,\qquad |i-j|=1,
\end{equation*}
$$F_i^2F_j-(q+q^{-1})F_iF_jF_i+F_jF_i^2=0,\qquad |i-j|=1,$$
$$E_iE_j-E_jE_I=F_iF_j-F_jF_i=0,\qquad |i-j|\ne 1.$$

The comultiplication $\Delta$, the antipode $S$, and the counit
$\varepsilon$ are determined by
\begin{equation*}
\Delta(E_i)=E_i \otimes 1+K_i \otimes E_i,\quad \Delta(F_i)=F_i
\otimes K_i^{-1}+1 \otimes F_i,\quad \Delta(K_i)=K_i \otimes K_i,
\end{equation*}
\begin{equation*}
S(E_i)=-K_i^{-1}E_i,\qquad S(F_i)=-F_iK_i,\qquad S(K_i)=K_i^{-1},
\end{equation*}
$$\varepsilon(E_i)=\varepsilon(F_i)=0,\qquad \varepsilon(K_i)=1.$$

The Hopf algebra $\mathbb{C}[SL_{N}]_q$ can be defined by the generators $t_{ij}$,
$i,j=1,...,N$, and the following relations
\begin{flalign*}
& t_{\alpha a}t_{\beta b}-qt_{\beta b}t_{\alpha a}=0, & a=b \quad \& \quad
\alpha<\beta,& \quad \text{or}\quad a<b \quad \& \quad \alpha=\beta,
\\ & t_{\alpha a}t_{\beta b}-t_{\beta b}t_{\alpha a}=0,& \alpha<\beta \quad
\&\quad a>b,& 
\\ & t_{\alpha a}t_{\beta b}-t_{\beta b}t_{\alpha a}-(q-q^{-1})t_{\beta a}
t_{\alpha b}=0,& \alpha<\beta \quad \& \quad a<b, & 
\\ & \det \nolimits_q \mathbf{t}=1.
\end{flalign*}
Here $\det_q \mathbf{t}$ is a q-determinant of the matrix
$\mathbf{t}=(t_{ij})_{i,j=1,\ldots,N}$:
\begin{equation*}\label{qdet}
\det \nolimits_q\mathbf{t}=\sum_{s \in
S_N}(-q)^{l(s)}t_{1s(1)}t_{2s(2)}\ldots t_{Ns(N)},
\end{equation*}
with $l(s)=\mathrm{card}\{(i,j)|\;i<j \quad \&\quad s(i)>s(j) \}$. The
comultiplication $\Delta$, the counit $\varepsilon$, and the antipode $S$
in $\mathbb{C}[SL_{N}]_q$ are defined as follows:
$$
\Delta(t_{ij})=\sum_kt_{ik}\otimes t_{kj},\qquad
\varepsilon(t_{ij})=\delta_{ij},\qquad S(t_{ij})=(-q)^{i-j}\det \nolimits_q
\mathbf{t}_{ji},
$$
with $\mathbf{t}_{ji}$ being the matrix derived from $\mathbf{t}$ by
discarding its $j$-th row and $i$-th column. This algebra is Hopf dual to $U_q \mathfrak{sl}_N$ and has the standard structure of $U_q^{\rm op}\mathfrak{sl}_N \otimes U_q \mathfrak{sl}_N$-module, where the superscript $\rm{'op'}$ means that the algebra $U_q^{\rm op}\mathfrak{sl}_N$ has the opposite multiplication (see \cite{Kl-Sch}). Moreover, the multiplication in $\mathbb{C}[SL_{N}]_q$ is a morphism of $U_q^{\rm op}\mathfrak{sl}_N \otimes U_q \mathfrak{sl}_N$-modules, i.e. $\mathbb{C}[SL_{N}]_q$ is a $U_q^{\rm op}\mathfrak{sl}_N \otimes U_q \mathfrak{sl}_N$-module algebra.

Recall that $U_q \mathfrak{sl}_N$ possesses a so called non-compact involution $*$ defined on generators as follows:
\begin{equation*}
(K_j^{\pm 1})^*=K_j^{\pm 1},\quad E_j^*=
\begin{cases}
K_jF_j,& j \ne n
\\ -K_jF_j,& j=n
\end{cases},\quad F_j^*=
\begin{cases}
E_jK_j^{-1},& j \ne n
\\ -E_jK_j^{-1},& j=n
\end{cases},\quad 
j=1,\ldots,N-1.
\end{equation*}
The corresponding $*$-Hopf algebra is denoted by $U_q \mathfrak{su}_{n,m}$ and may be treated as a quantum analog of the noncompact group $SU(n,m)$. The algebra  $\mathbb C[SL_N]_q$ possesses a compatible involution 
\begin{equation*}
t_{ij}^*= \mathrm{sign}\left((i-m-1/2)(n-j+1/2)\right)(-q)^{j-i}\det
\nolimits_q \mathbf{t}_{ij},
\end{equation*}
where  q-minors of $\mathbf{t}$ are defined in a natural way:
$$
t_{IJ}^{\wedge k}\stackrel{\mathrm{def}}{=}\sum_{s \in
S_k}(-q)^{l(s)}t_{i_1j_{s(1)}}\cdot t_{i_2j_{s(2)}}\cdots t_{i_kj_{s(k)}},
$$
with $I=\{(i_1,i_2,\dots,i_k)|\;1 \le i_1<i_2<\dots<i_k \le N \}$,
$J=\{(j_1,j_2,\dots,j_k)|\;1 \le j_1<j_2<\dots<j_k \le N \}$.
By means of that, $(\mathbb C[SL_N]_q,*)$ becomes a $U_q \mathfrak{su}_{n,m}$-module algebra, i.e. involutions in $U_q \mathfrak{sl}_N$ and $\mathbb C[SL_N]_q$ satisfy the following compatibility condition $(\xi f)^*=(S(\xi))^*f^*$ for every $\xi \in U_q \mathfrak{sl}_N, f \in \mathbb C[SL_N]_q$.

Let us also recall here basic notation in $q$-calculus. So,

$$(a;q)_k=\begin{cases}
(1-a)\cdot \ldots \cdot (1-aq^{k-1}), & k>0
\\ 1, & k=0.
\end{cases}$$
Eigenfunctions of our $q$-difference operators appear first as basic hypergeometric functions, so let us recall that a basic hypergeometric function is defined for $a_1,a_2,\ldots,a_r;b_1,b_2,\ldots,b_s\in\mathbb{C}$ as the infinite series
\begin{equation*}
\sideset{_r}{_s}{\mathop{\Phi}}\left[\begin{array}{c}a_1,a_2,\ldots,a_r\\
b_1,b_2,\ldots,b_s\end{array};q,z\right] =\sum_{k=0}^\infty
\frac{(a_1;q)_k(a_2;q)_k\cdots(a_r;q)_k}{(q;q)_k(b_1;q)_k\cdots(b_s;q)_k}
\left[(-1)^kq^{\frac{k(k-1)}2}\right]^{1+s-r}z^k.
\end{equation*}

We will also use the Jackson integral. This integral with the $q^{-2}$-base is defined by the formula $$\int_1^\infty f(x)d_{q^{-2}}x=\sum_{a=0}^{\infty} f(q^{-2a})q^{-2a}.$$

\section{\boldmath Algebras of functions on the quantum $\mathscr{H}_{n,m}$. Invariant integral}\label{sec_2}

Algebras of functions on a hyperbolic space were introduced in \cite{BS}. Here we just recall the notations and main objects along the way. 

We start with the well known \cite{RTF} $q$-analog of pseudo-Hermitian spaces. Let $\operatorname{Pol}(\widehat{\mathscr{H}}_{n,m})_q$ be the unital $*$-algebra with generators
$t_1,t_{2},\ldots,t_N$ and the commutation relations as follows:
\begin{equation*}
\begin{aligned}
t_it_j &= qt_jt_i,\qquad i<j,
\\ t_it_j^* &= qt_j^*t_i,\qquad i\ne j,
\\ t_it_i^* &= t_i^*t_i+(q^{-2}-1)\sum_{k=i+1}^Nt_kt_k^*,\qquad i>n,
\\ t_it_i^* &= t_i^*t_i+(q^{-2}-1)\sum_{k=i+1}^nt_kt_k^*-
(q^{-2}-1)\sum_{k=n+1}^Nt_kt_k^*,\qquad i\le n.
\end{aligned}
\end{equation*}
Obviously,
$$c=-\sum_{j=1}^nt_jt_j^*+\sum_{j=n+1}^Nt_jt_j^*$$
is a central element and not a zero divisor in $\operatorname{Pol}(\widehat{\mathscr{H}}_{n,m})_q$. Thus the localization
$\operatorname{Pol}(\widehat{\mathscr{H}}_{n,m})_{q,c}$ with respect to the multiplicative system
$c^{\mathbb{N}}$ is well defined.

The map $J:t_j \mapsto t_{1j}, j=1,\ldots,N,$ uniquely extends to an algebra homomorphism $J: \operatorname{Pol}(\widehat{\mathscr{H}}_{n,m})_{q,c} \rightarrow \mathbb C[SL_N]_q$. Note that $J(c)=\det_q \mathbf{t}=1$ since $J(c)$ is just the $q$-Laplace expansion of $\det_q \mathbf{t}$ along the first row.

The structure of $U_{q}\mathfrak{su}_{n,m}$-module algebra is transfered to $\operatorname{Pol}(\widehat{\mathscr{H}}_{n,m})_{q,c}$ from
the corresponding $U_{q}\mathfrak{su}_{n,m}$-module structure in $\mathbb{C}[SL_N]_q$ via the embedding $J$. 
 
Let $I_\phi$ be a $*$-automorphism of $\operatorname{Pol}(\widehat{\mathscr{H}}_{n,m})_{q,c}$ defined as follows: 
\begin{equation}\label{I_phi}
I_\phi(t_j)=e^{i\phi}t_j,\qquad j=1,2\ldots,N.
\end{equation}
Then
$$
\operatorname{Pol}(\mathscr{H}_{n,m})_q=\left\{\left.f\in
\operatorname{Pol}(\widehat{\mathscr{H}}_{n,m})_{q,c} \right|\: I_\phi(f)=f\right\}.
$$
is a $*$-algebra and will be called the algebra of regular functions on the quantum hyperbolic space. One can also distinguish this algebra as an algebra of $U_q \mathfrak{s(gl}_1 \times \mathfrak{gl}_{N-1})$-invariants in $\mathbb C[SL_N]_q$ w.r.t. the left action $\mathcal{L}$. 
\begin{equation}\label{left_invariance}
\operatorname{Pol}(\mathscr{H}_{n,m})_q=\left\{\left.f\in
\mathbb C[SL_N]_q \right|\: \mathcal{L}(\xi)f=\epsilon(\xi)f, \xi \in U_q \mathfrak{s(gl}_1 \times \mathfrak{gl}_{N-1}) \right\}.
\end{equation}

\medskip

Now we would like to recall the construction of an invariant integral on the quantum hyperbolic space. These results are also obtained in \cite{BS} and we give here just necessary definitions and statements. 

Recall that in the classical case there exists a positive invariant integral on $\mathscr{H}_{n,m}$, namely, one can integrate smooth functions with compact support on $\mathscr{H}_{n,m}$. So we need first to construct a $q$-analog for the corresponding space of functions. The construction uses a faithfull $*$-representaion of $\operatorname{Pol}(\mathscr{H}_{n,m})_q$ and we will start with it. 

The space $\mathscr{H}$ is a linear span of its basis
$\{e(i_1,i_2,\ldots,i_{N-1})|\:i_1,\ldots,i_n\in-\mathbb{Z}_+;\;
i_{n+1},\ldots,i_{N-1}\in\mathbb{N}\}$. The representation $T$ of 
$\operatorname{Pol}(\widehat{\mathscr{H}}_{m,n})_q$ is defined by
\begin{equation*} 
\begin{aligned}
T(t_j)e(i_1,\ldots,i_{N-1}) &=
q^{\sum\limits_{k=1}^{j-1}i_k}\cdot\left(q^{2(i_j-1)}-1\right)^{1/2}
e(i_1,\ldots,i_j-1,\ldots,i_{N-1}),
\\ T(t_j^*)e(i_1,\ldots,i_{N-1}) &=
q^{\sum\limits_{k=1}^{j-1}i_k}\cdot\left(q^{2i_j}-1\right)^{1/2}
e(i_1,\ldots,i_j+1,\ldots,i_{N-1}),
\end{aligned}
\end{equation*}
for $j\le n$,
\begin{equation*} 
\begin{aligned}
T(t_j)e(i_1,\ldots,i_{N-1}) &=
q^{\sum\limits_{k=1}^{j-1}i_k}\cdot\left(1-q^{2(i_j-1)}\right)^{1/2}
e(i_1,\ldots,i_j-1,\ldots,i_{N-1}),
\\ T(t_j^*)e(i_1,\ldots,i_{N-1}) &=
q^{\sum\limits_{k=1}^{j-1}i_k}\cdot\left(1-q^{2i_j}\right)^{1/2}
e(i_1,\ldots,i_j+1,\ldots,i_{N-1}),
\end{aligned}
\end{equation*}
for $n<j<N$, and, finally,
\begin{equation*} 
\begin{aligned}
T(t_N)e(i_1,\ldots,i_{N-1}) &=
q^{\sum\limits_{k=1}^{N-1}i_k}e(i_1,\ldots,i_{N-1}),
\\ T(t_N^*)e(i_1,\ldots,i_{N-1}) &=
q^{\sum\limits_{k=1}^{N-1}i_k}e(i_1,\ldots,i_{N-1}).
\end{aligned}
\end{equation*}
Obviously, one can restrict this representation to $\operatorname{Pol}(\mathscr{H}_{n,m})_q$.

The next proposition was proved in \cite{BS}.
\begin{proposition}
There exists a scalar product in $\mathscr{H}$ such that $T$ is a faithful $*$-representation of $\operatorname{Pol}(\mathscr{H}_{n,m})_q$ in the pre-Hilbert space $\mathscr{H}$.
\end{proposition}

Define the elements $\{x_j\}_{j=1,\ldots,N} \subset \operatorname{Pol}(\mathscr{H}_{n,m})_q$ as follows:
\begin{equation*}\label{x_j}
x_j\overset{\mathrm{def}}{=}
\begin{cases}
\sum\limits_{k=j}^Nt_kt_k^*, & j>n,
\\ -\sum\limits_{k=j}^nt_kt_k^*+\sum\limits_{k=n+1}^Nt_kt_k^*, & j\le n.
\end{cases}
\end{equation*}
Obviously, $x_1=1$, $x_ix_j=x_jx_i$, $i,j=1,\ldots,N$.

The vectors $e(i_1,\ldots,i_{N-1})$ are joint eigenvectors of the operators
$T(x_j)$, $j=1,2,\ldots,N$:
\begin{equation*}\label{Tx_jev}
\begin{aligned}
&T(x_1)=I,&
\\ &T(x_j)e(i_1,\ldots,i_{N-1})=
q^{2\sum\limits_{k=1}^{j-1}i_k}e(i_1,\ldots,i_{N-1}).
\end{aligned}
\end{equation*}

The joint spectrum of the pairwise commuting operators $T(x_j)$,
$j=1,2,\ldots,N$, is
\begin{multline*}\label{jsh}
\mathfrak{M}=\left\{(x_1,\ldots,x_N)\in\mathbb{R}^N\right|
\\ \left. x_i/x_j\in q^{2\mathbb{Z}}\;\&\;1=x_1\le x_2\le\ldots\le
x_{n+1}>x_{n+2}>\ldots>x_N>0\right\}.
\end{multline*}

\medskip

One can endow $\operatorname{Pol}(\widehat{\mathscr{H}}_{n,m})_q$ with the weakest topology such that the matrix elements are continuous. The completion of $\operatorname{Pol}(\widehat{\mathscr{H}}_{n,m})_q$ w.r.t. this topology will be considered as the space of generalized functions on the quantum $\widehat{\mathscr{H}}_{n,m}$ and
denoted by $\mathscr{D}(\widehat{\mathscr{H}}_{n,m})_q'$. Naturally, one can extend $T$ to a representation of $\mathscr{D}(\widehat{\mathscr{H}}_{n,m})_q'$ by continuity. Now one can identify $\mathscr{D}(\widehat{\mathscr{H}}_{n,m})_q'$ with the space of formal series
\begin{equation*}
f=\sum_{\fontsize{8}{1}(i_1,\ldots,i_N,j_1,\ldots,j_N):\;i_kj_k=0} t_1^{i_1}\ldots
t_n^{i_n}t_{n+1}^{*i_{n+1}}\ldots t_N^{*i_N}f_{IJ}(x_1,x_2,\ldots,x_N)t_N^{j_N}\ldots
t_{n+1}^{j_{n+1}}t_n^{*j_n}\ldots t_1^{*j_1},
\end{equation*}
where $f_{IJ}(x_1,x_2,\ldots,x_N)$ are functions on $\mathfrak{M}$. The topology on this space of formal series is the topology of pointwise convergence of the functions $f_{IJ}$.

Denote by $f_0$ the following function
\begin{equation*}\label{f_0(x_n+1)}
f_0=f_0(x_{n+1})=
\begin{cases}
1, & x_{n+1}=1,
\\ 0, & x_{n+1}\in q^{-2\mathbb{N}}.
\end{cases}
\end{equation*}
Introduce now a $*$-algebra $\operatorname{Fun}(\widehat{\mathscr{H}}_{n,m})_q \subset
\mathscr{D}(\widehat{\mathscr{H}}_{n,m})_q'$ generated by
$\operatorname{Pol}(\widehat{\mathscr{H}}_{n,m})_q$ and $f_0$. One can extend the $*$-automorphism $I_{\phi}$ (see \eqref{I_phi}) to a $*$-automorphism of the algebra $\operatorname{Fun}(\widehat{\mathscr{H}}_{n,m})_q$ by requiring $I_\phi(f_0)=f_0$. Then one can consider 
$$
\operatorname{Fun}(\mathscr{H}_{n,m})_q=\left\{\left.f\in
\operatorname{Fun}(\widehat{\mathscr{H}}_{n,m})_q \right|\: I_\phi(f)=f\right\}.
$$
Obviously, there exists a unique extension of the $*$-representation $T$ to a $*$-representation of  $\operatorname{Fun}(\mathscr{H}_{n,m})_q$ such that $T(f_0)$ is the orthogonal projection of $\mathscr{H}$ onto the submanifold corresponding to the first spectral point of $x_{n+1}$.

Let $\mathscr{D}(\mathscr{H}_{n,m})_{q,\mathfrak{k}}$ be the two-sided ideal
of $\operatorname{Fun}(\mathscr{H}_{n,m})_q$ generated by $f_0$. We call this
ideal the algebra of finite functions on the quantum hyperbolic space. It is
a quantum analog for the algebra of $K$-finite smooth functions on $\mathscr{H}_{n,m}$ with compact support.

\begin{remark}
It was explained in \cite{BS} that $\mathscr{D}(\mathscr{H}_{n,m})_{q,\mathfrak{k}}$ contains only such functions $\phi(x_{n+1})$, $x_{n+1} \in q^{-2\mathbb Z_+}$, that have finitely many nonzero values. We will also call such functions finite functions in variable $x=x_{n+1}$.
\end{remark}

Now we recall an explicit formula for a positive invariant integral on the space of finite functions
$\mathscr{D}(\mathscr{H}_{n,m})_{q,\mathfrak{k}}$.
Let $\nu_q:\mathscr{D}(\mathscr{H}_{n,m})_{q,\mathfrak{k}} \to\mathbb{C}$ be a linear functional defined by
\begin{equation*}
\nu_q(f)=\operatorname{Tr}(T(f)\cdot Q)=\int\limits_{\mathscr{H}_{n,m}} fd\nu_q,
\end{equation*}
where $Q=\mathrm{const}_1T(x_2 \ldots x_N)$ with $\mathrm{const}_1=(q^{-2};q^{-2})_{m-1}.$
\begin{theorem}[\cite{BS}]
The functional $\nu_q$ is well defined, positive, and $U_{q}\mathfrak{su}_{n,m}$-invariant.
\end{theorem}
Note that the constant is fixed by the requirement $\nu_q(f_0)=1$.

This invariant integral allows us to introduce a scalar product in $\mathscr{D}(\mathscr{H}_{n,m})_{q,\mathfrak{k}}$ by the formula $(f,g)=\int\limits_{\mathscr{H}_{n,m}} g^*f d\nu_q$.

The following proposition describes the structure of $\mathscr{D}(\mathscr{H}_{n,m})_{q,\mathfrak{k}}$ as a $U_q \mathfrak{k}$-module.
\begin{proposition}[\cite{BK1}]\label{Kdecomp_H}
As a $U_q \mathfrak{k}$-module, $\mathscr{D}(\mathscr{H}_{n,m})_{q,\mathfrak{k}}= \bigoplus\limits_{k,l,k',l' \in \mathbb Z_+, k+l'=l+k'} V_{(k,l,k',l')},$ where each isotypic component $V_{(k,l,k',l')}$ is a sum of infinite number of modules $L^{(n)}(k\varpi_1+l\varpi_{n-1}) \boxtimes L^{(m)}(l'\varpi_1+k'\varpi_{m-1})$ with $k+l'=l+k'$.
\end{proposition}

One can easily show that any highest weight vector from the isotypic component $V_{(k,l,k',l')}$ has the form $t_1^k(t_N^*)^{k'}\phi(x)t_{n+1}^{l'}(t_n^*)^l,$ where $\phi(x)$ is a finite function in $x$.

Now we would like to obtain restrictions of the scalar product $(\cdot,\cdot)$ to isotypic components $V_{(k,l,k',l')}$. By explicit computations, one may obtain the following result
\begin{lemma}\label{scalprod_highvect} For all $k,l,k',l' \in \mathbb Z_+$, $k+l'=l+k',$ one has
\begin{multline*}
(t_1^k(t_N^*)^{k'}\phi(x)t_{n+1}^{l'}(t_n^*)^l,t_1^k(t_N^*)^{k'}\psi(x)t_{n+1}^{l'}(t_n^*)^l)=C(k,l,k',l') \cdot \int_1^\infty \overline{\psi(x)}\phi(x)\rho_{k+l,k'+l'}(x)d_{q^{-2}}x,
\end{multline*}
where
\begin{equation*}
\rho_{a,b}(x)=x^{b+m-1}(q^{-2}x;q^{-2})_{a+n-1},
\end{equation*}
and
\begin{multline*} C(k,l,k',l')=(-1)^{k+l}q^{m(m-1)+2(k'+l')(m-1)-2lm}(q^{-2};q^{-2})_{m-1} \cdot \\ \frac{(q^2;q^2)_{k'}(q^2;q^2)_{l'}}{(q^2;q^2)_{k'+l'+m-1}}  \frac{(q^{-2};q^{-2})_{k}(q^{-2};q^{-2})_{l}}{(q^{-2};q^{-2})_{k+l+n-1}}.
\end{multline*}
\end{lemma}

The proof of this lemma is given in the Appendix.
\begin{corollary}
\begin{multline*}||t_1^k(t_N^*)^{k'}f_s(x)t_{n+1}^{l'}(t_n^*)^l ||^2=(-1)^{k+l}q^{m(m-1)-2lm+2(m-1)(k'+l')}(q^{-2};q^{-2})_{m-1} \\ \frac{(q^2;q^2)_{k'}(q^2;q^2)_{l'}}{(q^2;q^2)_{k'+l'+m-1}}\frac{(q^{-2};q^{-2})_{k}(q^{-2};q^{-2})_{l}}{(q^{-2};q^{-2})_{k+l+n-1}} q^{-2s(k'+l'+m)}(q^{-2s-2};q^{-2})_{k+l+n-1}.
\end{multline*}
\end{corollary}

\section{The Laplace-Beltrami operator and its restrictions to isotypic components}\label{sec_6}

In this section we introduce a $U_q \mathfrak{sl}_N$-invariant operator
$\Box_q$ in the space of finite functions $\mathscr{D}(\mathscr{H}_{n,m})_{q,
\mathfrak{k}}$. This operator will be considered as a quantum analog for the
invariant Laplace-Beltrami operator on the complex hyperbolic space
$\mathscr{H}_{n,m}$. Also we compute more explicitly its action on each isotypic component $V_{(k,l,k',l')}$. The restriction is closely related to a second order $q$-difference
operator in $x$.

First we should recall the definition of an operator $\overline{\partial}$ and a Hermitian pairing on 1-forms, as it was done in \cite{BK}.

The operator $\partial$ comes from the left action $\mathcal{L}$ of $U_q \mathfrak{sl}_N^{op}$ in $\operatorname{Pol}(\mathscr{H}_{n,m})_{q,c}$ embedded in $\mathbb C[SL_N]_q$, namely, $\partial=\mathcal{L}(K_1^{1/2}F_1)$. Then one can define $\overline{\partial}f=(\partial f^*)^*$. Let $\Omega^{(0,1)}(\mathscr{H}_{n,m})_q$ be the corresponding space of 1-forms, i.e. a module generated by $\operatorname{Pol}(\mathscr{H}_{n,m})_q$ and $\overline{\partial}t_{1i}$ inside $\mathbb C[SL_N]_q$.
A Hermitian pairing $(\cdot,\cdot)$ on $\Omega^{(0,1)}(\mathscr{H}_{n,m})_q$ is defined via the projection $P: \mathbb C[SL_N]_q \rightarrow \operatorname{Pol}(\mathscr{H}_{n,m})_q$ parallel to all non-zero $U_q \mathfrak{s(gl}_1 \times \mathfrak{gl}_{N-1})$-isotypic components of the representation $\mathcal{L}$ (recall that $\operatorname{Pol}(\mathscr{H}_{n,m})_q$ is the trivial isotypic component in $\mathbb C[SL_N]_q$ with respect to the left action of $U_q \mathfrak{s(gl}_1 \times \mathfrak{gl}_{N-1})$, see \eqref{left_invariance}). Evidently, the Hermitian pairing on (0,1)-forms satisfies the property $(\xi \omega_1,\omega_2)=(\omega_1,\xi^* \omega_2)$ for all $\xi \in U_q \mathfrak{sl}_N$, i.e. it is $U_q \mathfrak{sl}_N$-invariant. We set 
\begin{equation*}
(\theta_1,\theta_2)=P(\theta_2^*\theta_1), \qquad \theta_1,\theta_2 \in \Omega^{(0,1)}(\mathscr{H}_{n,m})_q.
\end{equation*}

Now we define $\Box_q$ by the formula
\begin{equation*}
\int_{\mathscr{H}_{n,m}}f_2^* (\Box_q f_1) d\nu_q=\int _{\mathscr{H}_{n,m}} (\bar \partial f_1, \bar
\partial f_2) d\nu_q, \qquad f_1,f_2 \in \mathscr{D}(\mathscr{H}_{n,m})_{q,\mathfrak{k}}.
\end{equation*}

\begin{proposition}
$\Box_q$ is a self-adjoint $U_q \mathfrak{sl}_N$-invariant operator.
\end{proposition}
{\bf Proof.} Let us check that $\Box_q$ is a $U_q \mathfrak{sl}_N$-invariant operator. Recall that the right and left $U_q \mathfrak{sl}_N$-actions in $\mathbb C[SL_N]_q$ commute with each other, thus $\bar \partial$ commute with the (right) $U_q \mathfrak{sl}_N$-action. Then we have for all $f_1,f_2 \in \mathscr{D}(\mathscr{H}_{n,m})_{q,\mathfrak{k}}$ and $\xi \in U_q \mathfrak{sl}_N$ 
\begin{equation*}
\int_{\mathscr{H}_{n,m}}f_2^* (\Box_q \xi f_1) d\nu_q= \int _{\mathscr{H}_{n,m}} (\bar \partial \xi f_1, \bar \partial f_2) d\nu_q= \\ \int _{\mathscr{H}_{n,m}} (\xi \bar \partial f_1, \bar \partial f_2) d\nu_q .
\end{equation*}
Due to the invariance of the Hermitian pairing $(\cdot,\cdot)$, we get 
\begin{multline*}
\int_{\mathscr{H}_{n,m}}f_2^* (\Box_q \xi f_1) d\nu_q= \int _{\mathscr{H}_{n,m}} (\bar \partial f_1, \xi^* \bar \partial f_2) d\nu_q =  \int _{\mathscr{H}_{n,m}} (\bar \partial f_1, \bar \partial \xi^* f_2) d\nu_q = \\ \int_{\mathscr{H}_{n,m}}(\xi^* f_2)^* (\Box_q f_1) d\nu_q = \int_{\mathscr{H}_{n,m}} (S(\xi)^{-1} f_2^*) (\Box_q f_1) d\nu_q = \int_{\mathscr{H}_{n,m}}f_2^* (\xi \Box_q f_1) d\nu_q.
\end{multline*}
The latter equality is due to the fact that $\nu_q$ is a $U_q \mathfrak{sl}_N$-invariant integral, hence it has an 'integration by parts property'
\begin{equation*}
\int_{\mathscr{H}_{n,m}} (\xi g_1) g_2 d\nu_q=\int_{\mathscr{H}_{n,m}} g_1 (S(\xi)g_2) d\nu_q , \quad g_1,g_2 \in \mathscr{D}(\mathscr{H}_{n,m})_{q,\mathfrak{k}}, \xi \in U_q \mathfrak{sl}_N.
\end{equation*}

Now we can verify that $\Box_q$ is a symmetric operator in $\mathscr{D}(\mathscr{H}_{n,m})_{q,\mathfrak{k}}$. Indeed,
\begin{multline*}
(\Box_q f_1,f_2)=\int_{\mathscr{H}_{n,m}}f_2^* (\Box_q f_1) d\nu_q= \int _{\mathscr{H}_{n,m}} (\bar \partial f_1, \bar \partial f_2) d\nu_q = \overline{\int _{\mathscr{H}_{n,m}} (\bar \partial f_2, \bar \partial f_1) d\nu_q}= \\ \overline{\int _{\mathscr{H}_{n,m}} f_1^* \Box_q f_2 d\nu_q}=\overline{(\Box_q f_2,f_1)}=(f_1, \Box_q f_2).
\end{multline*}
Since $\Box_q$ is defined on the whole space $\mathscr{D}(\mathscr{H}_{n,m})_{q,\mathfrak{k}}$, it is self-adjoint. \hfill $\blacksquare$

\medskip

Now we will study restrictions of $\Box_q$ to isotypic components $V_{(k,l,k',l')}$ of $\mathscr{D}(\mathscr{H}_{n,m})_{q,\mathfrak{k}}$. Due to the previous proposition, the restriction of $\Box_q$ to $V_{(k,l,k',l')}$ is completely defined by its action on the highest weight vectors $t_1^k(t_N^*)^{k'}\phi(x)t_{n+1}^{l'}(t_n^*)^l$.

Let us fix notation for $q$-difference operators:
\begin{equation*}
B^-: f(x) \mapsto \frac{f(q^{-2}x)-f(x)}{q^{-2}x-x}, \qquad B^+: f(x)
\mapsto \frac{f(q^2x)-f(x)}{q^2x-x}.
\end{equation*}

\begin{proposition}\label{Laplace_on_highestvector}
One has 
\begin{equation}\label{Box_hw}
\Box_q t_1^k(t_N^*)^{k'} \phi(x)t_{n+1}^{l'}(t_n^*)^l = t_1^k(t_N^*)^{k'} (A^{(k,l,k',l')}\phi(x))t_{n+1}^{l'}(t_n^*)^l ,
\end{equation} where
\begin{multline}\label{ALambda1}
A^{(k,l,k',l')} \phi (x)=q^{-1-2k-2l'}\frac{(1-q^{2(k+l')})(1-q^{2(N-1+k+l')})}{(1-q^2)(1-q^{2(N-1)})} \phi(x)- \\ q^{-1-2k'}\frac{1}{(1-q^2)(1-q^{2(N-1)})\rho_{k+l,k'+l'}(x)} B^+\bigg(\rho_{k+l,k'+l'}(x)x(q^{2(n+k)}-xq^{-2l})B^-\phi(x)\bigg).
\end{multline}
\end{proposition}

In order to expand the action of $\Box_q$ to the whole isotypic component, we simply need to apply elements from $U_q \mathfrak{k}$ to both sides of \eqref{Box_hw}. Hence all the information about the action of $\Box_q$ is stored in the introduced operators $A^{(k,l,k',l')}$.
The proof of this proposition is made by explicit calculations and will be given in the Appendix.

We can rewrite $A^{(k,l,k',l')}$ as a three terms $q$-difference operator in the following way
\begin{multline*}
A^{(k,l,k',l')} \psi(x)=\frac{q}{(1-q^2)(1-q^{2(N-1)})x}\cdot \big[ q^{-2k'-2l}(x-q^{2(n+k+l)})\psi(q^{-2}x) + \\ q^{2(N-2+k+l')}(x-1)\psi(q^2x)+ +q^{2(k-k'+n)}(1+q^{2(l'+k'+m-1)})\psi(x) -x(1+q^{2(N-1)})\psi(x)\big].
\end{multline*}

One can note that this operator depends only on the parameters $\Lambda=k+l$ and $\Lambda'=k'+l'$, since $k'+l=k+l'=\frac{k+l+k'+l'}{2}=\frac{\Lambda+\Lambda'}{2}$. Then 
\begin{multline*}\label{A_q_dif}
A^{\Lambda,\Lambda'} \psi(x)=  \frac{q}{(1-q^2)(1-q^{2(N-1)})x}\cdot \big[q^{-\Lambda-\Lambda'}(x-q^{2(n+\Lambda)})\psi(q^{-2}x) + q^{2N-2+\Lambda+\Lambda'}(x-1)\psi(q^2x)+\\+q^{2n+\Lambda-\Lambda'}(1+q^{2(m-1+\Lambda')})\psi(x) -x(1+q^{2(N-1)})\psi(x) \big].
\end{multline*}

\medskip

Let us consider the space of finite functions (functions that take only finitely many nonzero values at points from $q^{-2\mathbb Z_+}$) in one variable $x=x_{n+1}$ on $q^{-2\mathbb Z_+}$ with the inner product
\begin{equation*}
(\phi,\psi)_{\Lambda,\Lambda'}=\int_1^\infty \overline{\psi(x)}\phi(x) d\nu_q^{\Lambda,\Lambda'},
\end{equation*}
where $d\nu_q^{\Lambda,\Lambda'}= \frac{1}{(q^{-2};q^{-2})_{\Lambda+n-1}} \rho_{\Lambda,\Lambda'}(x) d_{q^{-2}}x$.
The scalar multiplier is chozen in such way that $(f_0,f_0)_{\Lambda,\Lambda'}=1$. 
We will complete this space of functions w.r.t. the corresponding norm and denote the resulting Hilbert space $L^2(d\nu_q^{\Lambda,\Lambda'})$. 

\begin{proposition}
$A^{\Lambda,\Lambda'}$ is a bounded self adjoint operator in $L^2(d\nu_q^{\Lambda,\Lambda'})$.
\end{proposition}
{\bf Proof.} Let us consider the operator $\frac{1}{\rho_{\Lambda,\Lambda'}(x)} B^+\bigg(\rho_{\Lambda,\Lambda'}(x)x(q^{2(n+k)}-xq^{-2l})B^-\phi(x)\bigg)$. It differs from the operator $x^{-(N-2+\Lambda+\Lambda')} B^+\bigg(x^{N+\Lambda+\Lambda'}B^-\phi(x)\bigg)$ by a compact operator. Thn one can check that the operator $x^{-a}B^+(x^{a+2}B^-\phi(x)),$ $a \in \mathbb N,$ is bounded, thus $A^{\Lambda,\Lambda'}$ is  bounded as well.

Let us check that $A^{\Lambda,\Lambda'}$ is a symmetric operator. Consider the finite functions on $q^{-2\mathbb Z_+}$
\[
f_j(x)=\begin{cases} 1, & x=q^{-2j},\\ 0, & \text{otherwise}, \end{cases} \qquad j \in \mathbb Z_+.
\]
These functions form an orthogonal system in $L^2(d\nu_q^{\Lambda,\Lambda'})$, so it is enough to verify symmetry on the basis vectors $f_j$.
One has
\begin{multline}\label{three_term_rec}
A^{\Lambda,\Lambda'}f_j= \frac{q}{(1-q^2)(1-q^{2(N-1)})x}\cdot \bigg( f_{j+1}(x-1)q^{2(N-1)+\Lambda+\Lambda'}+f_{j-1}(x-q^{2(n+\Lambda)})q^{-\Lambda-\Lambda'}\\ +f_j
(q^{2n+\Lambda-\Lambda'}+q^{2(N-1)+\Lambda+\Lambda'}-x(q^{2(N-1)}+1)) \bigg), \qquad j \in \mathbb Z_+.
\end{multline}
Thus $\int_1^\infty  \Big(A^{\Lambda,\Lambda'}f_j\Big) f_k d\nu^{\Lambda,\Lambda'}_q$ is nonzero only if $j=k$ or $j=k\pm 1$. Let us check the case $j=k-1$, the other two are also simple. For $j=k-1$ there is just one nonvanishing term from \eqref{three_term_rec} and
\begin{multline*}
\int_1^\infty  \Big(A^{\Lambda,\Lambda'}f_j \Big) f_k d\nu^{\Lambda,\Lambda'}_q= q^{2N-1+\Lambda+\Lambda'}\frac{q^{-2k}-1}{(1-q^2)(1-q^{2(N-1)})q^{-2k}}q^{-2k(\Lambda'+m)}(q^{-2-2k};q^{-2})_{\Lambda+n-1}= \\-q^{2N-1+\Lambda+\Lambda'}q^{-2k(\Lambda'+m-1)}\frac{(q^{-2k};q^{-2})_{\Lambda+n}}{(1-q^2)(1-q^{2(N-1)})}.
\end{multline*}
On the other side, for $k=j+1$ one has
\begin{multline*}
\int_1^\infty  f_j \Big(A^{\Lambda,\Lambda'} f_k\Big) d\nu^{\Lambda,\Lambda'}_q= q^{1-\Lambda-\Lambda'}\frac{q^{2-2k}-q^{2(n+\Lambda)}}{(1-q^2)(1-q^{2(N-1)})q^{2-2k}}q^{-2(k-1)(\Lambda'+m)}(q^{-2k};q^{-2})_{\Lambda+n-1}= \\ -q^{1-\Lambda-\Lambda'}q^{2(n+\Lambda)}q^{2(\Lambda'+m-1)}q^{-2k(\Lambda'+m-1)}\frac{(q^{-2k};q^{-2})_{\Lambda+n}}{(1-q^2)(1-q^{2(N-1)})},
\end{multline*}
so the scalar products are equal and $A^{\Lambda,\Lambda'}$ is a symmetric operator. $\hfill \blacksquare$

\section{Eigenfunctions of $A^{\Lambda,\Lambda'}$ and Al-Salam-Chihara polynomials}

\begin{lemma}\label{c1.5.2} The function
$$
\Phi^{\Lambda,\Lambda'}_l(x) = {_{3}\Phi_{2}}\left(
\begin{array}{c}
x,q^{\Lambda+\Lambda'-2l},q^{\Lambda+\Lambda'+2(l+N-1)}\\ q^{2(n+\Lambda)},0
\end{array}
;q^2,q^2\right)
$$
in $\mathscr{D}(\mathscr{H}_{n,m})_q^\prime$ is a generalized eigenfunction
for $A^{\Lambda,\Lambda'}$:
$$A^{\Lambda,\Lambda'} \Phi^{\Lambda,\Lambda'}_l=\lambda(l)\Phi^{\Lambda,\Lambda'}_l$$
with the eigenvalue
$$
\lambda(l) = - q\frac{(1-q^{-2l})(1-q^{2(l+N-1)})}{(1-q^2)(1-q^{2(N-1)})}.
$$
\end{lemma}

Recall the definition of the Al-Salam-Chihara polynomials, following \cite{koekoek}:
\begin{equation*}
Q_k(z;a,b|q)=\frac{(ab;q)_k}{a^k}\ {_{3}\Phi_{2}}\left(
\begin{array}{c}
q^{-k},a e^{i\theta},ae^{-i\theta}\\ ab,0
\end{array}
;q,q\right), \qquad z=\cos\theta.
\end{equation*}

We have the following well known orthogonality relations for the Al-Salam-Chihara polynomials:
\begin{multline}\label{rel_orth_ASC}
\frac{1}{2\pi}\int_{-1}^1 Q_i(z;a,b|q)Q_j(z;a,b|q)\frac{w(z)}{\sqrt{1-z^2}}\,dz+\\
\sum_{1<aq^k<a}w_k Q_i(z_k;a,b|q)Q_j(z_k;a,b|q)=\frac{\delta_{ij}}{(q^{i+1},abq^i;q)_\infty},
\end{multline}
where
\begin{equation*}
w(z)=\frac{h(z,1)h(z,-1)h(z,q^{1/2})h(z,-q^{1/2})}{h(z,a)h(z,b)}, \qquad
h(z,a)=(ae^{i\theta},ae^{-i\theta};q)_{\infty}, \,z=\cos\theta,
\end{equation*}
$z_k=\frac{aq^k+aq^{-k}}{2}$, and
\begin{equation*}
w_k=\frac{(a^{-2};q)_\infty}{(q,ab,b/a;q)_\infty}\frac{(1-a^2q^{2k})(a^2,ab;q)_k}
{(1-a^2)(q,qa/b;q)_k}q^{-k^2}\left(\frac{1}{a^3b}\right)^k.
\end{equation*}

An important fact about this family of polynomials is that there is a recurrence relation. Namely,
$$2zQ_k=Q_{k+1}+(a+b)q^{2k}Q_k+(1-q^{2k})(1-abq^{2k-2})Q_{k-1}.$$
We will employ this relation later on.

As one can see, the eigenfunctions $\Phi^{\Lambda,\Lambda'}_l(x)$ are connected with the Al-Salam-Chihara polynomials:
$$
\Phi^{\Lambda,\Lambda'}_l(q^{-2k}) = \frac{q^{k(N-1+\Lambda+\Lambda')}}{(q^{2(n+\Lambda)};q^2)_k}Q_k(z;q^{2n-N+1+\Lambda-\Lambda'},q^{N-1+\Lambda+\Lambda'}|q^2), \qquad
e^{i\theta}=q^{2l+N-1},\, z=\cos\theta.
$$
This relation gives us a proof of Lemma \ref{c1.5.2}.
Indeed, for every $k \in \mathbb Z_+$ one has
\begin{multline*}
A^{\Lambda,\Lambda'} \Phi^{\Lambda,\Lambda'}_l(q^{-2k})=q^{1-\Lambda-\Lambda'}\frac{q^{-2k}-q^{2(n+\Lambda)}}{(1-q^2)(1-q^{2(N-1)})q^{-2k}}\Phi^{\Lambda,\Lambda'}_l(q^{-2k-2}) +\\ q^{2N-1+\Lambda+\Lambda'}\frac{q^{-2k}-1}{(1-q^2)(1-q^{2(N-1)})q^{-2k}}\Phi^{\Lambda,\Lambda'}_l(q^{-2k+2})+\\+q^{2n+1+\Lambda-\Lambda'}\frac{1+q^{2(m-1+\Lambda')}}{(1-q^2)(1-q^{2(N-1)})q^{-2k}}\Phi^{\Lambda,\Lambda'}_l(q^{-2k}) -q\frac{(1+q^{2(N-1)})}{(1-q^2)(1-q^{2(N-1)})}\Phi^{\Lambda,\Lambda'}_l(q^{-2k}).
\end{multline*}
Thus, after substituting the expression for the Al-Salam-Chihara polynomials, one gets
\begin{multline*}
A^{\Lambda,\Lambda'} Q_k(z;q^{2n-N+1+\Lambda-\Lambda'},q^{N-1+\Lambda+\Lambda'}|q^2) =\\ \frac{q^{N}}{(1-q^2)(1-q^{2(N-1)})} Q_{k+1}(z;q^{2n-N+1+\Lambda-\Lambda'},q^{N-1+\Lambda+\Lambda'}|q^2) + \\ \frac{q^{N}(1-q^{2k})(1-q^{2(n+\Lambda+k-1)})}{(1-q^2)(1-q^{2(N-1)})}Q_{k-1}(z;q^{2n-N+1+\Lambda-\Lambda'},q^{N-1+\Lambda+\Lambda'}|q^2)+\\+ \frac{q^{2k}(q^{2N-1+\Lambda+\Lambda'}+q^{2n+1+\Lambda-\Lambda'})-q(1+q^{2(N-1)})}{(1-q^2)(1-q^{2(N-1)})}Q_k(z;q^{2n-N+1+\Lambda-\Lambda'},q^{N-1+\Lambda+\Lambda'}|q^2).
\end{multline*}
One can notice that the r.h.s. of this formula is similar to the r.h.s. of recurrence relation for the Al-Salam-Chihara polynomials. Thus one may continue calculations as follows
 
\begin{multline*}
A^{\Lambda,\Lambda'} Q_k(z;q^{2n-N+1+\Lambda-\Lambda'},q^{N-1+\Lambda+\Lambda'}|q^2) = \\ 2z \cdot \frac{q^{N}}{(1-q^2)(1-q^{2(N-1)})} Q_k(z;q^{2n-N+1+\Lambda-\Lambda'},q^{N-1+\Lambda+\Lambda'}|q^2) - \\ \frac{q(1+q^{2(N-1)})}{(1-q^2)(1-q^{2(N-1)})}Q_k(z;q^{2n-N+1+\Lambda-\Lambda'},q^{N-1+\Lambda+\Lambda'}|q^2).
\end{multline*}
Recall that $z=\cos \theta=\frac{q^{2l+N-1}+q^{-(2l+N-1)}}{2}$, so
the eigenvalue 
\begin{multline*}
\lambda(l)=(q^{2l+N-1}+q^{-(2l+N-1)})\cdot \frac{q^{N}}{(1-q^2)(1-q^{2(N-1)})}-\frac{q(1+q^{2(N-1)})}{(1-q^2)(1-q^{2(N-1)})}= \\ -q \frac{(1-q^{-2l})(1-q^{2(l+N-1)})}{(1-q^2)(1-q^{2(N-1)})}. \hfill \blacksquare
\end{multline*}

Let us denote the orthogonal measure for Al-Salam-Chihara polynomials by $d\sigma^{\Lambda,\Lambda'}$.  

Recall that $q=e^{-h/2}$. One may note that if $z=\cos \theta$ satistfies $|z|\leq 1$, then $2l+N-1 \in i\mathbb R$, i.e. $2l=-N+1+i\nu,$ $\nu \in \mathbb R$.
Moreover, one can make a change of variables and rewrite the continuous part of the orthogonal measure as follows
\begin{multline*}
\frac{h}{2\pi}\int_{0}^{\pi/h} Q_i(z;a,b|q)Q_j(z;a,b|q)w(z)\,d\nu+\\
\sum_{1<aq^k<a}w_k Q_i(z_k;a,b|q)Q_j(z_k;a,b|q)=\frac{\delta_{ij}}{(q^{i+1},abq^i;q)_\infty},
\end{multline*}
where $z=\frac{q^{i\nu}+q^{-i\nu}}{2}$ and 
\begin{equation*}
w(z)=\Big|\frac{(e^{2i\theta};q^2)_\infty}{(ae^{i\theta},be^{i\theta};q^2)_\infty}\Big|^2=\Big|\frac{(q^{2i\nu};q^2)_\infty}{(q^{i\nu+N-1+\Lambda+\Lambda'},q^{i\nu +n-m+1+\Lambda-\Lambda'};q^2)_\infty}\Big|^2.
\end{equation*}
Let us denote 
$$c_{\Lambda,\Lambda'}(a)=\frac{(q^{a+N-1+\Lambda+\Lambda'},q^{a +n-m+1+\Lambda-\Lambda'};q^2)_\infty}{(q^{2a};q^2)_\infty},
$$
so that $w(z)=\Big|\frac{1}{c_{\Lambda,\Lambda'}(i\nu)}\Big|^2$ for $z=\frac{q^{i\nu}+q^{-i\nu}}{2}$. Thus we get a $q$-analog for the Harish-Chandra c-function, which appear in the Plancherel measure for complex hyperbolic spaces (see \cite{Faraut}).

\section{A spectral theorem for $A^{\Lambda,\Lambda'}$}\label{sec_9}

In this section we obtain a spectral theorem for $A^{\Lambda,\Lambda'}$. The obtained spectral measures contribute to the Plancherel measure of the whole Laplace-Beltrami operator. As in the classical case (see
\cite[p.429-432]{Faraut}), the support of the Plancherel measure consists of continuous and discrete parts. This structure reflects the decomposition of a quasi-regular representation of $U_q \mathfrak{su}_{n,m}$ on quantum complex hyperbolic space into series of unitary $U_q \mathfrak{su}_{n,m}$-modules.

\begin{theorem}
The bounded self-adjoint linear operator $A^{\Lambda,\Lambda'}$ is unitary equivalent to the operator of multiplication
by independent variable in the Hilbert space $L^2(d\sigma^{\Lambda,\Lambda'})$. The unitary equivalence is given by the operator
\begin{align*}
U^{\Lambda,\Lambda'}:\, L^2(d\nu_q^{\Lambda,\Lambda'}) &\rightarrow L^2(d\sigma^{\Lambda,\Lambda'}),
\\U^{\Lambda,\Lambda'}:\,f(x)&\mapsto\hat{f}(\lambda)=\int\limits_1^\infty f(x)\Phi_l(x) d \nu_q^{\Lambda,\Lambda'},
\end{align*}
where $\lambda=-q\frac{(1-q^{-2l})(1-q^{2(l+N-1)})}{(1-q^2)(1-q^{2(N-1)})}$.
\end{theorem}

{\bf Proof.}  By standard arguments \cite{AG},
the bounded self-adjoint linear operator $A^{\Lambda,\Lambda'}$ is unitary equivalent to the multiplication operator
$f(\lambda) \mapsto \lambda f(\lambda)$ in the Hilbert space $L^2(d\mu^{\Lambda,\Lambda'}(\lambda))$ of square integrable
functions with respect to a certain measure $d\mu^{\Lambda,\Lambda'}(\lambda)$ with compact support in $\mathbb R$. Let us find
explicitly the corresponding measure and the operator of unitary equivalence $U^{\Lambda,\Lambda'}$. One can fix the unitary
equivalence operator by the condition $U^{\Lambda,\Lambda'}f_0=1$.

By easy calculations, $||f_j||^2=q^{-2j(N-1+\Lambda+\Lambda')}\frac{(q^{2j+2};q^2)_{\Lambda+n-1}}{(q^2;q^2)_{\Lambda+n-1}}$. Thus the finite functions $e_j=q^{j(N-1+\Lambda+\Lambda')}\sqrt{\frac{(q^{2};q^2)_{\Lambda+n-1}}{(q^{2j+2};q^2)_{\Lambda+n-1}}}f_j$ form an
orthonormal system in $L^2(d\nu_q^{\Lambda,\Lambda'})$. Due to \eqref{three_term_rec}, $A^{\Lambda,\Lambda'}$ acts on them by the formula
\begin{gather*}
A^{\Lambda,\Lambda'} e_j=\frac{q^{N}}{(1-q^2)(1-q^{2(N-1)})}(
e_{j+1}\sqrt{(1-q^{2j+2})(1-q^{2j+2n+2\Lambda})} + \\ e_{j-1}\sqrt{(1-q^{2j})(1-q^{2j+2n+2\Lambda-2})} + \\e_j
(q^{2j+(N-1)+\Lambda+\Lambda'}+q^{2j+2n-(N-1)+\Lambda-\Lambda'}-q^{N-1}-q^{1-N})), \qquad j \in \mathbb Z_+.
\end{gather*}
Thus $P^{\Lambda,\Lambda'}_j=U^{\Lambda,\Lambda'}e_j \in L^2(d\mu^{\Lambda,\Lambda'}(\lambda))$, $j\in\mathbb{Z}_+$, form an orthonormal system of polynomials
\[
\int P^{\Lambda,\Lambda'}_i(\lambda)P^{\Lambda,\Lambda'}_j(\lambda)d\mu^{\Lambda,\Lambda'}(\lambda)=\delta_{ij},\qquad i,j\in\mathbb{Z}_+,
\]
and one has
\begin{multline}\label{P_rec}
\lambda P^{\Lambda,\Lambda'}_j(\lambda)= \frac{q^{N}}{(1-q^2)(1-q^{2(N-1)})}
(P^{\Lambda,\Lambda'}_{j+1}(\lambda)\sqrt{(1-q^{2j+2})(1-q^{2j+2n+2\Lambda})} +\\ P^{\Lambda,\Lambda'}_{j-1} (\lambda) \sqrt{(1-q^{2j})(1-q^{2j+2n+2\Lambda-2})}+\\
P^{\Lambda,\Lambda'}_j(\lambda) (q^{2j+(N-1)+\Lambda+\Lambda'}+q^{2j+2n-(N-1)+\Lambda-\Lambda'}-q^{N-1}-q^{1-N})), \qquad j \in \mathbb Z_+,
\end{multline}
\begin{equation}\label{init_P0}
P^{\Lambda,\Lambda'}_0(\lambda)=U^{\Lambda,\Lambda'}e_0=U^{\Lambda,\Lambda'}f_0=1.
\end{equation}
(we naturally suppose that $P^{\Lambda,\Lambda'}_{-1}=0$).

The orthogonal polynomials $P^{\Lambda,\Lambda'}_j(\lambda)$, $j\in\mathbb{Z}_+$, are determined by \eqref{P_rec} and
\eqref{init_P0}. Let us compare them with the corresponding recurrence relations and the
initial data for the Al-Salam-Chihara polynomials $Q_j(z;q^{2n-(N-1)+\Lambda-\Lambda'},q^{N-1+\Lambda+\Lambda'}|q^2)$. One can observe that
\[
||Q_j||^2=\frac{1}{(q^{2j+2},q^{2j+2n+2\Lambda};q^2)_\infty}=\frac{(q^2,q^{2n+2\Lambda};q^2)_j}{(q^2,q^{2n+2\Lambda};q^2)_\infty}.
\]

The polynomials
\[
\widetilde{Q}_j\overset{\mathrm{def}}{=}\sqrt{\frac{1}{(q^2,q^{2n+2\Lambda};q^2)_j}} \,Q_j
\]
satisfy the following recurrence relations:
\begin{gather*}
z\widetilde{Q}_j=\dfrac12\sqrt{(1-q^{2(j+1)})(1-q^{2j+2n+2\Lambda})}\widetilde{Q}_{j+1}+
\dfrac12\sqrt{(1-q^{2j})(1-q^{2j+2n+2\Lambda-2})}\widetilde{Q}_{j-1}+\\
\dfrac12q^{2j}(q^{2n-(N-1)+\Lambda-\Lambda'}+q^{N-1+\Lambda+\Lambda'})\widetilde{Q}_j, \qquad j \in \mathbb Z_+.
\end{gather*}

Thus we obtain that $P^{\Lambda,\Lambda'}_j$ and $\widetilde{Q}_j$ are related by the change of variable
\[
\lambda=\frac{2q^{N}}{(1-q^2)(1-q^{2(N-1)})}(2z-q^{N-1}-q^{1-N}).
\]
So,
\[
P^{\Lambda,\Lambda'}_j(\lambda)=\sqrt{\frac{1}{(q^2,q^{2n+2\Lambda};q^2)_j}} \,Q_j(z;q^{2n-N+1+\Lambda-\Lambda'},q^{N-1+\Lambda+\Lambda'}|q^2),
\]
where $z=\frac{q^{2l+N-1}+q^{-(2l+N-1)}}{2}$ and
\[\lambda=-\frac{q}{(1-q^2)(1-q^{2(N-1)})}(1-q^{2l})(1-q^{-2(l+N-1)}).\]

\begin{multline*}
U^{\Lambda,\Lambda'}f_j=q^{-j(N-1+\Lambda+\Lambda')}\sqrt{\frac{(q^{2j+2};q^2)_{\Lambda+n-1}}{(q^{2};q^2)_{\Lambda+n-1}}}\sqrt{\frac{1}{(q^2,q^{2n+2\Lambda};q^2)_j}}\,
Q_j(z;q^{2n-N+1+\Lambda-\Lambda'},q^{N-1+\Lambda+\Lambda'}|q^2) \\ =q^{-j(N-1+\Lambda+\Lambda')}\frac{1}{(q^2;q^2)_j}\,Q_j(z;q^{2n-N+1+\Lambda-\Lambda'},q^{N-1+\Lambda+\Lambda'}|q^2).
\end{multline*}
On the other hand,
\[
\int\limits_1^\infty
f_j(x)\Phi^{\Lambda,\Lambda'}_{l}(x)\rho_{\Lambda,\Lambda'}(x)d_{q^{-2}}x=\frac{q^{-j(N-1+\Lambda+\Lambda')}}{(q^2;q^2)_j}\,Q_j(z;q^{2n-N+1+\Lambda-\Lambda'},q^{N-1+\Lambda+\Lambda'}|q^2),
\]
where $z=\frac12(q^{2l+N-1}+q^{-(2l+N-1)})$.

Hence for every function $f(x)$ on $q^{-2\mathbb{Z}_+}$ with finite support one has
\[
U^{\Lambda,\Lambda'}f=\int\limits_1^\infty f(x)\Phi^{\Lambda,\Lambda'}_{l}(x)\rho_{\Lambda,\Lambda'}(x)d_{q^{-2}}x.
\]
Now the claim of the theorem follows from the orthogonality relations for the Al-Salam-Chihara polynomials
\eqref{rel_orth_ASC}. \hfill $\blacksquare$

\begin{remark}
\begin{equation*}
\hat{f_j}=q^{-2j(\Lambda+\Lambda'+N-1)}\frac{(q^{2j+2};q^2)_{\Lambda+n-1}}{(q^2;q^2)_{\Lambda+n-1}}\Phi^{\Lambda,\Lambda'}_l(q^{-2j}).
\end{equation*}
\end{remark}
\bigskip

\section{Appendix. Proof of technical results}

{\bf Proof of Lemma \ref{scalprod_highvect}.} In order to compute the scalar product of the two highest weight vectors $f_1,f_2$ from the same isotypic component $V_{(k,l,k',l')}$, first we need to find how the operator $T(f_2^*f_1)$ acts in the Fock space.

Let $f_1=t_1^k(t_N^*)^{k'}\phi(x)t_{n+1}^{l'}(t_n^*)^l$, $f_2=t_1^k(t_N^*)^{k'}\psi(x)t_{n+1}^{l'}(t_n^*)^l$. Then explicit calculations show that
\begin{multline*}
T(f_2^*f_1) e(a_1\ldots a_{N-1})= (-1)^{k+l} q^{2ll'+2k'l-2k'l'} (q^{2a_1-2};q^{-2})_k (q^{2a_n};q^2)_l q^{-2la_n} \\ (q^{2a_{n+1}-2};q^{-2})_{l'} q^{2k'\sum_{j=1}^{N-1} a_j} q^{2(l+l') \sum_{j=1}^n a_j} \overline{\psi(q^{2l+2\sum_{j=1}^n a_j})} \phi(q^{2l+2\sum_{j=1}^n a_j}) e(a_1 \ldots a_{N-1}).
\end{multline*}
Then 
\begin{multline*}
\mathrm{Tr}T(f_2^*f_1) Q= (-1)^{k+l} q^{2ll'+2k'l-2k'l'} \sum_{a_1,\ldots, a_{N-1}} (q^{2a_1-2};q^{-2})_k (q^{2a_n};q^2)_l q^{-2la_n} \\ (q^{2a_{n+1}-2};q^{-2})_{l'} q^{2k'\sum_{j=1}^{N-1} a_j} q^{2(l+l') \sum_{j=1}^n a_j} \overline{\psi(q^{2l+2\sum_{j=1}^n a_j})} \phi(q^{2l+2\sum_{j=1}^n a_j}) q^{2\sum_{j=1}^{N-1}(N-j)a_j}= 
\\ (-1)^{k+l} q^{2ll'+2k'l-2k'l'} \sum_{a_1,\ldots,a_n} (q^{2a_1-2};q^{-2})_k (q^{2a_n};q^2)_l q^{-2la_n} q^{2(l+l') \sum_{j=1}^n a_j}  \overline{\psi(q^{2l+2\sum_{j=1}^n a_j})} \phi(q^{2l+2\sum_{j=1}^n a_j}) \\ \sum_{a_{n+1},\ldots, a_{N-1}} (q^{2a_{n+1}-2};q^{-2})_{l'} q^{2\sum_{j=1}^{N-1}(N-j+k')a_j}.
\end{multline*}
Let us reorganise the first sum as the sum over the paramter $c=\sum_{j=1}^n a_j$ and the sum over $a_1,\ldots,a_n$ subject to the condition that their sum should be equal to some integer $c$. Note that the latter sum then will be finite, because all $a_1,\ldots,a_n$ are nonpositive integers. Then one can continue the calculations as follows
\begin{multline*}
\mathrm{Tr}T(f_2^*f_1) Q= (-1)^{k+l} q^{2ll'+2k'l-2k'l'} \sum_{c=-\infty}^0 q^{2c(l+l'+k'+m)}\overline{\psi(q^{2l+2c})}\phi(q^{2l+2c}) \\ \sum_{a_1,\ldots,a_n, \sum a_j=c} (q^{2a_1-2k};q^2)_k (q^{2a_n};q^2)_l q^{-2la_n} q^{2\sum_{j=1}^n (n-j)a_j} \sum_{a_{n+1},\ldots, a_{N-1}} (q^{2a_{n+1}-2};q^{-2})_{l'} q^{2\sum_{j=n+1}^{N-1}(N-j+k')a_j}.
\end{multline*}
Applying summation formulas from Lemma \ref{summation_lemma}, we get
\begin{multline*}
\mathrm{Tr}T(f_2^*f_1) Q= (-1)^{k+l} q^{2ll'+2k'l-2k'l'} \sum_{c=-\infty}^0 q^{2c(l+l'+k'+m)} \overline{\psi(q^{2l+2c})} \phi(q^{2l+2c}) \\ (q^{-2};q^{-2})_k (q^{-2};q^{-2})_l q^{-2lc} \frac{(q^{2(c+l-1)};q^{-2})_{k+l+n-1}}{(q^{-2};q^{-2})_{k+l+n-1}} \cdot q^{2k'l'}q^{m(m-1)+2(m-1)(k'+l')} \frac{(q^2;q^2)_{k'}(q^2;q^2)_{l'}}{(q^2;q^2)_{k'+l'+m-1}}=
\\ (-1)^{k+l} q^{2ll'+2k'l} q^{m(m-1)+2(m-1)(k'+l')} \frac{(q^{-2};q^{-2})_k (q^{-2};q^{-2})_l}{(q^{-2};q^{-2})_{k+l+n-1}}\frac{(q^2;q^2)_{k'}(q^2;q^2)_{l'}}{(q^2;q^2)_{k'+l'+m-1}}\cdot 
\\ \sum_{c=-\infty}^0 q^{2c(l'+k'+m)} \overline{\psi(q^{2l+2c})} \phi(q^{2l+2c}) (q^{2(c+l-1)};q^{-2})_{k+l+n-1}
\end{multline*}
One can notice that in the last sum if $c \geq 1-l$, then the corresponding term vanishes. So the sum is equal to
\begin{multline*}
\sum_{c=-\infty}^{-l} q^{2c(l'+k'+m)} \overline{\psi(q^{2l+2c})} \phi(q^{2l+2c}) (q^{2(c+l-1)};q^{-2})_{k+l+n-1}= \\
q^{-2l(k'+l'+m)}\sum_{c=-\infty}^0 q^{2c(l'+k'+m)} \overline{\psi(q^{2c})} \phi(q^{2c}) (q^{2(c-1)};q^{-2})_{k+l+n-1}= \\ q^{-2l(k'+l'+m)} \int_{1}^\infty \overline{\psi(x)} \phi(x) x^{k'+l'+m-1}(q^{-2}x;q^{-2})_{k+l+n-1}d_{q^{-2}}x.
\end{multline*}
 
Finally, we get 
\begin{multline*}
\mathrm{Tr}T(f_2^*f_1) Q= (-1)^{k+l} q^{m(m-1)+2(m-1)(k'+l')-2lm} \frac{(q^{-2};q^{-2})_k (q^{-2};q^{-2})_l}{(q^{-2};q^{-2})_{k+l+n-1}} \cdot \frac{(q^2;q^2)_{k'}(q^2;q^2)_{l'}}{(q^2;q^2)_{k'+l'+m-1}}\cdot 
\\ \int_{1}^\infty \overline{\psi(x)} \phi(x) x^{k'+l'+m-1}(q^{-2}x;q^{-2})_{k+l+n-1}d_{q^{-2}}x. \hfill \blacksquare
\end{multline*}
\begin{remark}
Although we claimed Lemma \ref{scalprod_highvect} only for highest weight vectors from $\mathscr{D}(\mathscr{H}_{n,m})_{q,\mathfrak{k}}$, i.e. there was a requirement that $k+l'=k'+l$, one can easily see that $$f_2^*f_1=t_n^{l}(t_{n+1}^*)^{l'}t_N^{k'}(t_1^*)^kt_1^k(t_N^*)^{k'}\phi(x)\psi(x)t_{n+1}^{l'}(t_n^*)^{l} \in \mathscr{D}(\mathscr{H}_{n,m})_{q,\mathfrak{k}}$$ for all integral values of $k,l,k',l'$. Thus the integral of such expressions is well defined, and of course the answer is given by the same formula.
\end{remark}

The next lemma contains auxilarly summation formulas which allows us explicitly compute invariant integrals in the Lemma \ref{scalprod_highvect}. 
\begin{lemma}\label{summation_lemma}
\begin{enumerate}
\item 
\begin{multline*}\sum \limits_{a_1+\ldots+a_n=-t} (q^{2a_1-2k};q^{2})_k(q^{2a_n};q^2)_l q^{-2la_n} q^{2\sum_{i=1}^n (n-i)a_i}= \\ (q^{-2};q^{-2})_k(q^{-2};q^{-2})_l q^{2lt} \frac{(q^{-2(t-l+1)};q^{-2})_{k+l+n-1}}{(q^{-2};q^{-2})_{k+l+n-1}}.
\end{multline*}
\item \begin{equation*} \sum \limits_{a_{n+1}, \ldots a_{N-1} \in \mathbb N} (q^{2a_{n+1}-2};q^{-2})_{l'} q^{2k'\sum_{i=n+1}^{N-1}a_i} q^{2\sum_{i=n+1}^{N-1}(N-i)a_i}= q^{(m-1)(2k'+2l'+m)}q^{2l'k'}\frac{(q^2;q^2)_{k'}(q^2;q^2)_{l'}}{(q^2;q^2)_{k'+l'+m-1}}.\end{equation*}
\end{enumerate}
\end{lemma}
{\bf Proof.} Let us prove the first summation formula. One can divide both sides of the required identity by $q^{2lt}(q^{-2};q^{-2})_k(q^{-2};q^{-2})_l$. Then certain $q$-binomial coefficients will appear; namely, one has
\begin{equation*}
\frac{(q^{2a_1-2k};q^{2})_{k}}{(q^{-2};q^{-2})_k}=\frac{(q^{-2};q^{-2})_{k-a_1}}{(q^{-2};q^{-2})_k(q^{-2};q^{-2})_{-a_1}} \qquad \frac{(q^{2a_n};q^{2})_l}{(q^{-2};q^{-2})_l}=\frac{(q^{-2};q^{-2})_{-a_n}}{(q^{-2};q^{-2})_l(q^{-2};q^{-2})_{-a_n-l}}
\end{equation*}
(recall that $a_1,a_n \leq 0$). Note that we can restrict ourselves to the summands with $a_n \leq -l$, otherwise $(q^{2a_n};q^{2})_l$ and the corresponding summand vanishes.
Recall the standard notation for the $q$-binomial coefficient with the $q^{-2}$-base
\begin{equation*}
\begin{bmatrix}
a \\ b
\end{bmatrix}_{q^{-2}}=\frac{(q^{-2};q^{-2})_{a}}{(q^{-2};q^{-2})_{b}(q^{-2};q^{-2})_{a-b}}.
\end{equation*}
Now one can rewrite the summation as follows
\begin{multline*}
\sum \limits_{a_1+\ldots+a_n=-t} \begin{bmatrix}
k-a_1 \\ k
\end{bmatrix}_{q^{-2}} \begin{bmatrix}
-a_n \\ l
\end{bmatrix}_{q^{-2}} q^{-2la_n} q^{2\sum_{i=1}^n (n-i)a_i}=\\ \sum \limits_{0 \geq a_1+a_n \geq -t} \begin{bmatrix}
k-a_1 \\ k
\end{bmatrix}_{q^{-2}} q^{2(n-1)a_1}\begin{bmatrix}
-a_n \\ l
\end{bmatrix}_{q^{-2}} q^{-2la_n}q^{-2lt} \sum_{a_2+\ldots+a_{n-1}=-t-a_1-a_n}q^{2\sum_{i=2}^{n-1} (n-i)a_i}.
\end{multline*}

The last sum was computed in \cite[Proposition ]{BK1}. The answer is again given by a $q$-binomial coefficient
\begin{equation*}
\sum_{a_2+\ldots+a_{n-1}=-t-a_1-a_n}q^{2\sum_{i=2}^{n-1} (n-i)a_i}=\begin{bmatrix}
n-3+t+a_1+a_n \\ n-3
\end{bmatrix}_{q^{-2}}.
\end{equation*}
Thus we need to compute the following sum
\begin{multline*}
\sum \limits_{0 \geq a_1+a_n \geq -t} \begin{bmatrix}
k-a_1 \\ k
\end{bmatrix}_{q^{-2}} q^{2(n-1)a_1}\begin{bmatrix}
-a_n \\ l
\end{bmatrix}_{q^{-2}} q^{-2la_n}q^{-2lt} \sum_{a_2+\ldots+a_{n-1}=-t-a_1-a_n}q^{2\sum_{i=2}^{n-1} (n-i)a_i} =\\ 
\sum \limits_{0 \geq a_1+a_n \geq -t} \begin{bmatrix}
k-a_1 \\ k
\end{bmatrix}_{q^{-2}} q^{2(n-1)a_1}\begin{bmatrix}
-a_n \\ l
\end{bmatrix}_{q^{-2}} q^{-2la_n} q^{-2lt}\begin{bmatrix}
n-3+t+a_1+a_n \\ n-3
\end{bmatrix}_{q^{-2}}=\\ \sum \limits_{0 \leq x+y \leq t} \begin{bmatrix}
k+x \\ k
\end{bmatrix}_{q^{-2}} q^{-2(n-1)x}\begin{bmatrix}
y \\ l
\end{bmatrix}_{q^{-2}} q^{2ly}q^{-2lt} \begin{bmatrix}
n-3+t-x-y \\ n-3
\end{bmatrix}_{q^{-2}}
\end{multline*}

Now this summation can be computed by applying twice the following formula
\begin{equation*}
\sum_{x+y=t}\begin{bmatrix}
k+x \\ k
\end{bmatrix}_{q^{-2}}\begin{bmatrix}
l+y \\ l
\end{bmatrix}_{q^{-2}} q^{-2x(l+1)} = \begin{bmatrix}
k+l+t+1 \\ k+l+1
\end{bmatrix}_{q^{-2}}.
\end{equation*}
This identity is proved by induction in $t$.

Let us prove the second identity. One can easily observe that the summation in variables $a_{n+1}, \ldots, a_{N-1}$ are independent. Thus, one can rewrite the whole sum as a double sum
\begin{multline*}
\sum \limits_{a_{n+1}, \ldots a_{N-1} \in \mathbb N} (q^{2a_{n+1}-2};q^{-2})_{l'} q^{2k'\sum_{i=n+1}^{N-1}a_i} q^{2\sum_{i=n+1}^{N-1}(N-i)a_i}= \\ \sum_{a_{n+1 }=1}^\infty (q^{2a_{n+1}-2};q^{-2})_{l'} q^{2(N-n-1+k')a_{n+1}}\cdot \sum \limits_{a_{n+2}, \ldots a_{N-1} \in \mathbb N}  q^{2k'\sum_{i=n+2}^{N-1}a_i} q^{2\sum_{i=n+2}^{N-1}(N-i)a_i}.
\end{multline*}
Now the second sum is just a multibasic geometric progression, and this gives us
\begin{multline*}
\sum \limits_{a_{n+1}, \ldots a_{N-1} \in \mathbb N} (q^{2a_{n+1}-2};q^{-2})_{l'} q^{2k'\sum_{i=n+1}^{N-1}a_i} q^{2\sum_{i=n+1}^{N-1}(N-i)a_i}= \\ \sum_{a=1}^\infty (q^{2a-2};q^{-2})_{l'} q^{2(m-1+k')a}\cdot \prod \limits_{j=1}^{m-2}  \frac{q^{2(j+k')}}{1-q^{2(j+k')}}= \frac{q^{(2k'+m-1)(m-2)}}{(q^{2k'+2};q^2)_{m-2}} \sum_{a=1}^\infty (q^{2a-2};q^{-2})_{l'} q^{2(m-1+k')a}.
\end{multline*}
Now one can verify that
\begin{equation*}
 \sum_{a=1}^\infty (q^{2a-2};q^{-2})_x q^{2ya}=q^{2y(x+1)}\frac{(q^2;q^2)_x}{(q^{2y};q^2)_{x+1}}.
\end{equation*}
 Substituting the last formula, one gets that
\begin{multline*}
\sum \limits_{a_{n+1}, \ldots a_{N-1} \in \mathbb N} (q^{2a_{n+1}-2};q^{-2})_{l'} q^{2k'\sum_{i=n+1}^{N-1}a_i} q^{2\sum_{i=n+1}^{N-1}(N-i)a_i}= \\  \frac{q^{(2k'+m-1)(m-2)}}{(q^{2k'+2};q^2)_{m-2}} \cdot  q^{2(k'+m-1)(l'+1)}\frac{(q^2;q^2)_{l'}}{(q^{2(m-1+k')};q^2)_{l'+1}}= \\ q^{m(m-1)}q^{2(m-1)(k'+l')+2k'l'}\frac{(q^2;q^2)_{k'}(q^2;q^2)_{l'}}{(q^2;q^2)_{k'+m-2}(q^{2(m-1+k')};q^2)_{l'+1}},
\end{multline*}
which coinsides with the required answer.
\bigskip

Let us recall the notation for the weight function $\rho_{\Lambda,\Lambda'}(x)=x^{k'+l'+m-1}(q^{-2}x;q^{-2})_{k+l+n-1},$ where $\Lambda=k+l$ and $\Lambda'=k'+l'$.

\medskip

{\bf Proof of Proposition \ref{Laplace_on_highestvector}} We need to prove the equality 
\begin{equation*}
\Box_q t_1^k(t_N^*)^{k'} \phi(x)t_{n+1}^{l'}(t_n^*)^l = t_1^k(t_N^*)^{k'} (A^{(k,l,k',l')}\phi(x))t_{n+1}^{l'}(t_n^*)^l ,
\end{equation*}
where the operator $A^{k,l,k',l'}$ is given by \eqref{ALambda1}. Since $\Box_q$ is a $U_q \mathfrak{sl}_N$-invariant operator, it is obvious that $\Box_q (t_1^k(t_N^*)^{k'} \phi(x)t_{n+1}^{l'}(t_n^*)^l) \in \mathscr{D}(\mathscr{H}_{n,m})_{q,\mathfrak{k}}$ is again a highest weight vector of the same isotypic component. Thus, the action can be expressed as some operator on finite function $\phi(x)$, and our goal is to find explicit formula for it.
Let us start with computing the action of $\bar{\partial}$ on $t_1^k(t_N^*)^{k'} \phi(x)t_{n+1}^{l'}(t_n^*)^l.$ Explicit calculations show that 
\begin{equation*}
\bar{\partial} \phi(x)=\bar{\partial}x \cdot B^-\phi(x), \qquad B^-\phi(x)=\frac{\phi(q^{-2}x)-\phi(x)}{q^{-2}x-x}.
\end{equation*}
Denote $y=\bar{\partial}x$ and $s=k+l'=l+k'$. Then one has 
\begin{multline*}
\bar{\partial}(t_1^k(t_N^*)^{k'} \phi(x)t_{n+1}^{l'}(t_n^*)^l)= q^{1/2}\frac{1-q^{2l}}{1-q^2}t_1^k(t_N^*)^{k'} \phi(x)t_{n+1}^{l'}(t_n^*)^{l-1}t_{2,n}^*+\\ q^{1/2}q^{l+k-s}t_1^k(t_N^*)^{k'}y B^- \phi(x)t_{n+1}^{l'}(t_n^*)^l+q^{1/2}q^{l+k-s}\frac{1-q^{2k'}}{1-q^2}t_1^k(t_N^*)^{k'-1}t_{2,N}^* \phi(x)t_{n+1}^{l'}(t_n^*)^l.
\end{multline*}
Now we can plug in this formula to the scalar product $(\bar \partial t_1^k(t_N^*)^{k'}\phi(x)t_{n+1}^{l'}(t_n^*)^l, \bar \partial t_1^k(t_N^*)^{k'}\psi(x)t_{n+1}^{l'}(t_n^*)^l)$ and obtain that
\begin{multline*}
(\bar \partial t_1^k(t_N^*)^{k'}\phi(x)t_{n+1}^{l'}(t_n^*)^l, \bar \partial t_1^k(t_N^*)^{k'}\psi(x)t_{n+1}^{l'}(t_n^*)^l)=\\
q \Big(\frac{1-q^{2l}}{1-q^2}\Big)^2P(t_{2n}t_n^{l-1}(t_{n+1}^*)^{l'}\overline{\psi(x)} t_N^{k'}(t_1^*)^k t_1^k(t_N^*)^{k'}\phi(x)t_{n+1}^{l'}(t_n^*)^{l-1}t_{2n}^*)+\\
q^{1+l+k-s}\frac{1-q^{2l}}{1-q^2}P(t_{2n}t_n^{l-1}(t_{n+1}^*)^{l'}\overline{\psi(x)}t_N^{k'}(t_1^*)^k t_1^k(t_N^*)^{k'}yB^-\phi(x)t_{n+1}^{l'}(t_n^*)^{l})+\\
q^{1+l+k-s}\frac{1-q^{2l}}{1-q^2}\frac{1-q^{2k'}}{1-q^2}P(t_{2n}t_n^{l-1}(t_{n+1}^*)^{l'}\overline{\psi(x)}t_N^{k'}(t_1^*)^k t_1^k(t_N^*)^{k'-1}t_{2N}^*\phi(x)t_{n+1}^{l'}(t_n^*)^{l})+\\
q^{1+l+k-s}\frac{1-q^{2l}}{1-q^2}P(t_n^{l}(t_{n+1}^*)^{l'}\overline{B^-\psi(x)}y^*t_N^{k'}(t_1^*)^k t_1^k(t_N^*)^{k'}\phi(x)t_{n+1}^{l'}(t_n^*)^{l-1}t_{2n}^*)+\\
q^{1+2(l+k-s)}P(t_n^{l}(t_{n+1}^*)^{l'}\overline{B^-\psi(x)}y^*t_N^{k'}(t_1^*)^k t_1^k(t_N^*)^{k'}yB^-\phi(x)t_{n+1}^{l'}(t_n^*)^{l})+\\
q^{1+2(l+k-s)}\frac{1-q^{2k'}}{1-q^2}P(t_n^{l}(t_{n+1}^*)^{l'}\overline{B^-\psi(x)}y^*t_N^{k'}(t_1^*)^k t_1^k(t_N^*)^{k'-1}t_{2N}^*\phi(x)t_{n+1}^{l'}(t_n^*)^{l})+\\
q^{1+l+k-s}\frac{1-q^{2k'}}{1-q^2}\frac{1-q^{2l}}{1-q^2}P(t_n^{l}(t_{n+1}^*)^{l'}\overline{\psi(x)}t_{2N}t_N^{k'-1}(t_1^*)^k t_1^k(t_N^*)^{k'}\phi(x)t_{n+1}^{l'}(t_n^*)^{l-1}t_{2n}^*)+\\
q^{1+2(l+k-s)}\frac{1-q^{2k'}}{1-q^2}P(t_n^{l}(t_{n+1}^*)^{l'}\overline{\psi(x)}t_{2N}t_N^{k'-1}(t_1^*)^k t_1^k(t_N^*)^{k'}yB^-\phi(x)t_{n+1}^{l'}(t_n^*)^{l})+\\
q^{1+2(l+k-s)}\Big(\frac{1-q^{2k'}}{1-q^2}\Big)^2P(t_n^{l}(t_{n+1}^*)^{l'}\overline{\psi(x)}t_{2N}t_N^{k'-1}(t_1^*)^k t_1^k(t_N^*)^{k'-1}t_{2N}^*\phi(x)t_{n+1}^{l'}(t_n^*)^{l}).
\end{multline*}
Now let us recall that $P$ is the projection  $P: \mathbb C[SL_N]_q \rightarrow \operatorname{Pol}(\mathscr{H}_{n,m})_q \subset \mathbb C[SL_N]_q$ parallel to all non-zero isotypic components of the representation $\mathcal{L}$. It was verified in \cite[Lemma 3]{BK1} that
\begin{equation*}
P(t_{2i}t_{2j}^*)=q^{-2}\frac{1-q^2}{1-q^{2(N-1)}}(\epsilon_{ij}-t_{1i}t_{1j^*}), \qquad \epsilon_{ij}=\begin{cases}
0, & i\neq j,
\\ q^{2(i-1)}, & i=j \geq n+1, 
\\ -q^{2(i-1)}, & i=j \leq n.
\end{cases}
\end{equation*}

Applying these formulas and commutational relations in $\mathbb C[SL_N]_q$, one gets that
\begin{multline*}
(\bar \partial t_1^k(t_N^*)^{k'}\phi(x)t_{n+1}^{l'}(t_n^*)^l, \bar \partial t_1^k(t_N^*)^{k'}\psi(x)t_{n+1}^{l'}(t_n^*)^l)= 
\\ q^{-1} \frac{1-q^2}{1-q^{2(N-1)}} \Big(-q^{2n-2l} \left(\frac{1-q^{2l}}{1-q^2} \right)^2 t_n^{l-1}(t_{n+1}^*)^{l'}t_N^{k'}(t_1^*)^kt_1^k(t_N^*)^{k'}\overline{\psi(x)}\phi(x)t_{n+1}^{l'}(t_n^*)^{l-1}+ 
\\ q^{2(N+k-2k')} \left(\frac{1-q^{2k'}}{1-q^2} \right)^2 t_n^l(t_{n+1}^*)^{l'}t_N^{k'-1}(t_1^*)^kt_1^k(t_N^*)^{k'-1}\overline{\psi(x)}\phi(x)t_{n+1}^{l'}(t_n^*)^l + 
\\t_n^l(t_{n+1}^*)^{l'}t_N^{k'}(t_1^*)^kt_1^k(t_N^*)^{k'} \{ q^{-2-4k'}x(q^{2(n+k)}-xq^{-4l})\overline{B^-\psi(x)}B^-\phi(x) 
\\ +\frac{q^{-4k'}}{1-q^2} (q^{2(n+k)}(1-q^{2k'})-xq^{-2l}(q^{-2l}-q^{2k'}))(\overline{B^-\psi(x)} \phi(x)+\overline{\psi(x)}B^-\phi(x)) - 
\\q^{2-4k'} \left(\frac{q^{-2l}-q^{2k'}}{1-q^2} \right)^2\overline{\psi(x)}\phi(x) \} t_{n+1}^{l'}(t_n^*)^l \Big)
\end{multline*}

Now we need to compute the integral of this expression over the quantum $\mathscr{H}_{n,m}$. To do that, we apply Lemma \ref{scalprod_highvect}. So for $f_1=t_1^k(t_N^*)^{k'}\phi(x)t_{n+1}^{l'}(t_n^*)^l$, $f_2=t_1^k(t_N^*)^{k'}\psi(x)t_{n+1}^{l'}(t_n^*)^l$ one has $\Box_q f_1=t_1^k(t_N^*)^{k'}A^{(k,l,k',l')}\phi(x)t_{n+1}^{l'}(t_n^*)^l$ and
\begin{equation*}
\int_{\mathscr{H}_{n,m}}f_2^*(\Box_q f_1)d\nu_q=C(k,l,k',l')\int_1^\infty \overline{\psi(x)}A^{(k,l,k',l')}\phi(x) \rho_{\Lambda,\Lambda'}(x)d_{q^{-2}}x.
\end{equation*}
Let us now integrate $(\bar{\partial}f_1,\bar{\partial}f_2)$. Then one has
\begin{multline*}
\int_{\mathscr{H}_{n,m}}(\bar{\partial}f_1,\bar{\partial}f_2) d\nu_q= q^{-1} \frac{1-q^2}{1-q^{2(N-1)}} \cdot \\ \Big(-q^{2n-2l} \left(\frac{1-q^{2l}}{1-q^2} \right)^2 C(k,l-1,k',l') \int_1^\infty \overline{\psi(x)}\phi(x) \rho_{\Lambda-1,\Lambda'}(x) d_{q^{-2}}x+ 
\\ q^{2(N+k-2k')} \left(\frac{1-q^{2k'}}{1-q^2} \right)^2 C(k,l,k'-1,l') \int_1^\infty \overline{\psi(x)}\phi(x) \rho_{\Lambda,\Lambda'-1}(x)d_{q^{-2}}x+ 
\\ C(k,l,k',l') \int_1^\infty  A(\phi,\psi)(x) \rho_{\Lambda,\Lambda'}(x) d_{q^{-2}}x \Big),
\end{multline*}
where 
\begin{multline*}
A(\phi,\psi)(x)=q^{-2-4k'}x(q^{2(n+k)}-xq^{-4l})\overline{B^-\psi(x)}B^-\phi(x) 
+\frac{q^{-4k'}}{1-q^2} (q^{2(n+k)}(1-q^{2k'})- \\ xq^{-2l}(q^{-2l}-q^{2k'})) (\overline{B^-\psi(x)} \phi(x)+\overline{\psi(x)}B^-\phi(x)) - 
q^{2-4k'} \left(\frac{q^{-2l}-q^{2k'}}{1-q^2} \right)^2\overline{\psi(x)}\phi(x).
\end{multline*}
Further calculations show that 
\begin{multline*}
\int_{\mathscr{H}_{n,m}}(\bar{\partial}f_1,\bar{\partial}f_2) d\nu_q= q^{-1} \frac{C(k,l,k',l')}{(1-q^2)(1-q^{2(N-1)})} \cdot \\ \bigg(q^{2n-2l} (1-q^{2l})^2 \frac{1-q^{-2(k+l+n-1)}}{1-q^{-2l}}\int_1^\infty \overline{\psi(x)}\phi(x) \rho_{\Lambda-1,\Lambda'}(x) d_{q^{-2}}x+ 
\\ q^{2(N+k-2k')} (1-q^{2k'})^2 q^{-2(m-1)} \frac{1-q^{2(k'+l'+m-1)}}{1-q^{2k'}} \int_1^\infty \overline{\phi(x)}\phi(x) \rho_{\Lambda,\Lambda'-1}(x) d_{q^{-2}}x+ 
\\ (1-q^2)^2 \int_1^\infty  A(\phi,\psi)(x) \rho_{\Lambda-1,\Lambda'}(x) d_{q^{-2}}x \bigg)=
\\ q^{-1} \frac{C(k,l,k',l')}{(1-q^2)(1-q^{2(N-1)})} \cdot \\ \bigg(q^{-2(k+l-1)} (1-q^{2l})(1-q^{2(k+l+n-1)}) \int_1^\infty \overline{\psi(x)}\frac{\phi(x)}{1-q^{-2(k+l+n-1)}x}\rho_{\Lambda,\Lambda'}(x) d_{q^{-2}}x+ 
\\ q^{2(n+1+k-2k')} (1-q^{2k'})(1-q^{2(k'+l'+m-1)}) \int_1^\infty \overline{\psi(x)}\frac{\phi(x)}{x} \rho_{\Lambda,\Lambda'}(x) d_{q^{-2}}x+ 
\\ (1-q^2)^2 \int_1^\infty  A(\phi,\psi)(x)
\rho_{\Lambda,\Lambda'}(x) d_{q^{-2}}x \bigg).
\end{multline*}
Using $q$-analogs of the partial integration formulas one can prove that operators $B^-$ and $-q^2 B^+$ are formally dual. Namely, for every
functions $u(x), v(x)$ with finite support on $q^{2\mathbb{Z}}$ one has 
$$
\int\limits_0^\infty u(x) \frac{v(q^{-2}x)-v(x)}{(q^{-2}-1)x} d_{q^{-2}} x =
-q^2 \int\limits_0^\infty \frac{u(x)-u(q^2x)}{(1-q^2)x} v(x) d_{q^{-2}} x,
$$
where
$$
\int\limits_0^\infty f(x) d_{q^{-2}} x = (q^{-2} - 1)
\sum\limits_{k=-\infty}^\infty f(q^{-2k}) q^{-2k}.
$$
Thus one can transfer 'all activity' in the last integral from $\psi(x)$ to $\phi(x)$ and get that 
\begin{multline*}
\int_1^\infty  A(\phi,\psi)(x) \rho_{\Lambda,\Lambda'}(x) d_{q^{-2}}x = - q^{-4k'}\int_1^\infty \overline{\psi(x)}\cdot B^+\Big(x(q^{2(n+k)}-xq^{-4l})B^-\phi(x) \rho_{\Lambda,\Lambda'}(x)\Big) d_{q^{-2}} x\\ - \frac{q^{2-4k'}}{1-q^2}\int_1^\infty \overline{\psi(x)} B^+((q^{2(n+k)}(1-q^{2k'})-xq^{-2l}(q^{-2l}-q^{2k'}))\phi(x) \rho_{\Lambda,\Lambda'}(x)) d_{q^{-2}} x \\+\frac{q^{-4k'}}{1-q^2}\int_{1}^{\infty} \overline{\psi(x)} (q^{2(n+k)}(1-q^{2k'})-xq^{-2l}(q^{-2l}-q^{2k'}))B^-\phi(x) \rho_{\Lambda,\Lambda'}(x)) d_{q^{-2}} x \\ -q^{2-4k'} \left(\frac{q^{-2l}-q^{2k'}}{1-q^2} \right)^2\int_1^\infty \overline{\psi(x)}\phi(x)\rho_{\Lambda,\Lambda'}(x) d_{q^{-2}} x.
\end{multline*}
Hence one has
\begin{multline*}
A^{(k,l,k',l')} \phi (x)= \frac{q^{-1}}{(1-q^{2(N-1)})(1-q^2)}\Big(q^{2(m-l-k+1)}\frac{(1-q^{2l})(1-q^{2(l+k+n-1)})}{(1-q^{-2(l+k+n-1)}x)}\phi(x) + \\ q^{2(n-k'-l'+l+1)}\frac{(1-q^{2k'})(1-q^{2(l'+k'+m-1)})}{x}\phi (x)
\\ - q^{2-4k'} (q^{-2l}-q^{2k'})^2 \phi (x) +q^{-4k'}(1-q^2) (q^{2(n+k)}(1-q^{2k'})-xq^{-2l}(q^{-2l}-q^{2k'}))B^-\phi(x) \\-\frac{q^{2-4k'}(1-q^2)}{\rho_{\Lambda,\Lambda'}(x)} B^+ \left((q^{2(n+k)}(1-q^{2k'})-xq^{-2l}(q^{-2l}-q^{2k'}))\phi(x)\rho_{\Lambda,\Lambda'}(x) \right) - \\ \frac{q^{-4k'}(1-q^2)^2}{\rho_{\Lambda,\Lambda'}(x)} B^+(x(q^{2(n+k)}-xq^{-4l})\rho_{\Lambda,\Lambda'}(x)B^-\phi(x))\Big). 
\end{multline*}
What remains is to rewrite this expression in a shorter form. One can verify that, indeed, the last expression is equal to 
\begin{multline*}
A^{(k,l,k',l')} \phi (x)=q^{-1-2k-2l'}\frac{(1-q^{2(k+l')})(1-q^{2(N-1+k+l')})}{(1-q^2)(1-q^{2(N-1)})} \phi(x)- \\ q^{-1-2k'}\frac{q^{2n+2k}}{(1-q^2)(1-q^{2(N-1)})\rho_{\Lambda,\Lambda'}(x)} B^+\bigg(\rho_{\Lambda,\Lambda'}(x)x(1-xq^{-2(n+k+l)})B^-\phi(x)\bigg).\hfill \blacksquare
\end{multline*}


\bigskip

\begin{thebibliography}{40}
\ifx\undefined\BibUrl\def\BibUrl#1{\url{#1}}\else\fi
\ifx\undefined\BibAnnote\long\def\BibAnnote#1{}\else\fi
\ifx\undefined\BibEmph\def\BibEmph#1{\emph{#1}}\else\fi

\bibitem{AG} N.I.~Akhiezer, I.M.~Glazman, {\it Theory of Linear Operators in Hilbert Spaces}, Dover Publ. (1993).

\bibitem{BK} O.~Bershtein, Ye.~Kolisnyk, {\it Plancherel measure for the
quantum matrix ball - I}, Math. Phys., Anal., Geom. {\bf 5}, No.4, (2009), 315-346.

\bibitem{BK1} O.~Bershtein, Ye.~Kolisnyk, {\it Harmonic analysis on quantum complex hyperbolic spaces}, SIGMA {\bf 7} (2011), 078, 19 p.

\bibitem{BS} O.~Bershtein, S.~Sinelshchikov, {\it Function theory on a q-analog of complex hyperbolic space}, J. of Geometry and Physics {\bf 62}, No.5, (2012), 1323-1337.
math.QA/1009.6063.

\bibitem{Dijk-Sh} G.~van Dijk, Yu.~Sharshov, The Plancherel formula
for line bundles on complex hyperbolic spaces, {\it J. Math. Pures Appl.} {\bf 79} (2000), No. 5, 451--473.

\bibitem{RTF} L.D.~Faddeev, N.Yu.~Reshetikhin, L.A.~Takhtadjan, {\it
Quantization of Lie groups and Lie algebras}, Algebra and Analysis {\bf 1} No 1 (1989), 178 -- 206.

\bibitem{Faraut} J.~Faraut, {\it Distributions sph\'eriques sur les espaces hyperboliques}, J. Math. pures et appl. {\bf 58} (1979), 369 -- 444.

\bibitem{GR} G.~Gaspar, M.~Rahman, Basic Hypergeometric Series, Cambridge Univ. Press, Cambridge, (1990).

\bibitem{I} I.C.-H.~Ip {\it Representation of the quantum plane, its quantum double, and harmonic analysis on $GL^+_q(2,R)$,}  Sel. Math. New Ser. {\bf 19} (2013), 987--1082, doi:10.1007/s00029-012-0112-4.

\bibitem{Jant} J.C.~Jantzen, {\it Lectures on Quantum Groups,} Providence, R.I.: American Mathematical Society, (1996).

\bibitem{Kl-Sch} A.~Klimyk, K.~Schm\"udgen, {\it Quantum Groups and Their Representations,} Springer, Berlin et al., (1997).

\bibitem{KS} E.~Koelink, J.V.~Stokman, {\it Fourier transforms on the quantum $SU(1; 1)$ group, with an appendix by M. Rahman,}
Publ. Res. Inst. Math. Sci. 37 (2001), 621--715, math.QA/9911163.

\bibitem{koekoek} R.~Koekoek, R.F.~Swarttouw, {\it The Askey-scheme of
hypergeometric orthogonal polynomials and its $q$-analogue}, Report 98-17,
Technical University Delft, Delft, 1998.

\bibitem{Molch1} V.~Molchanov, {\it Spherical functions on hyperboloids}, Math.Sb. {\bf 99} (1976), No.2, 139--161. (Engl.transl.: Math. USSR-Sb., {\bf 28} (1976), 119--139.)

\bibitem{Molch2} V.~Molchanov, {\it Harmonic analysis on homogeneous spaces}, Itogi nauki i tekhn., Sovr.probl.mat.
Fund.napr. {\bf 59}, VINITI (1990), 5--144. Engl.transl.: Encycl. Math. {\bf 59}, Springer, Berlin etc. (1995), 1--135.

\bibitem{Schl} H.~Schlichtkrull, {\it Eigenspaces of the Laplacian on hyperbolic spaces: composition series and integral transforms}, J. Func. Anal. {\bf 70} (1987), 194--219.

\bibitem{polmat} D.~Shklyarov, S.~Sinel'shchikov, L.~Vaksman, {\it
Fock representations and quantum matrices}, International J. Math. {\bf 15} No.9 (2004), 855--894.


\bibitem{SZ} D.~Shklyarov, G.~Zhang, {\it Berezin transform on the quantum unit ball}, J. of Math. Phys. {\bf 44} 4344 (2003); http://doi.org/10.1063/1.1593226.

\bibitem{U} K.~Ueno, {\it Spectral analysis for the Casimir operator on the quantum group $SU(1; 1)$}, Proc. Japan Acad. Ser. A Math. Sci. {\bf 66} (1990), no. 2, 42--44.

\bibitem{V} L.L.~Vaksman, {\it Quantum Bounded Symmetric Domains}, translated by O.~Bershtein and
S.~Sinel'shchikov, Translations of the Mathematical Monographs, vol.238, AMS, (2010).

\bibitem{KV} L.L.~Vaksman, L.I.~Korogodskii, {\it Spherical functions on the quantum group $SU(1; 1)$ and a q-analogue of the Mehler-Fock formula,} Funktsional. Anal. i Prilozhen. {\bf 25} (1991), no. 1, 60--62 (English transl.: Funct. Anal. Appl. {\bf 25} (1991), no. 1, 48--49).

\end{thebibliography}
\end{document}